\documentclass[12pt]{article}
\usepackage{amssymb,amsmath}
\begin{document}
\begin{center}
{\Large\bf Interior capacities of condensers with~infinitely many
plates in a locally compact space}\\
\bigskip
{\large\it Natalia Zorii}\\
\end{center}
\bigskip

{\small{\bf Abstract.} The study deals with the theory of interior
capacities of condensers in a locally compact space, a condenser
being treated here as a countable, locally finite collection of
arbitrary sets with the sign~$+1$ or~$-1$ prescribed such that the
closures of opposite-signed sets are mutually disjoint. We are
motivated by the known fact that, in the noncompact case, the main
minimum-problem of the theory is in general unsolvable, and this
occurs even under very natural assumptions (e.\,g., for the
Newtonian, Green, or Riesz kernels in~$\mathbb R^n$, $n\geqslant2$,
and closed condensers of finitely many plates). Therefore it was
particularly interesting to find statements of variational problems
dual to the main min\-imum-prob\-lem (and hence providing some new
equivalent definitions of the capacity), but now always solvable
(e.\,g., even for nonclosed, unbounded condensers of infinitely many
plates). For all positive definite kernels satisfying B.~Fuglede's
condition of consistency between the strong and the vague
($={}$weak$*$) topologies, problems with the desired properties are
posed and solved. Their solutions provide a natural generalization
of the well-known notion of interior capacitary distributions
associated with a set. We give a description of those solutions,
establish statements on their uniqueness and continuity, and point
out their characteristic properties.\bigskip

\bigskip

}\bigskip

{\bf\large 1. Introduction}\medskip

The present work is devoted to further development of the theory of
interior capacities of condensers in a locally compact space. A
condenser will be treated here as a countable, locally finite
collection of arbitrary (noncompact or even nonclosed) sets with the
sign~$+1$ or~$-1$ prescribed such that the closures of
opposite-signed sets are mutually disjoint. For a background of the
theory for condensers of finitely many plates we refer the reader
to~\cite{Z2}--\cite{Z7}; see also~\cite{O}, where the condensers
were additionally assumed to be compact.

The reader is expected to be familiar with the principal notions and
results of the theory of measures and integration on a locally
compact space; its exposition can be found in~\cite{B2,E2} (see
also~\cite{F1,Z3} for a brief survey).

The theory of interior capacities of condensers provides a natural
extension of the well-known theory of interior capacities of sets,
developed by~H.~Cartan~\cite{Car} and Vall\'{e}e-Poussin~\cite{VP}
for classical kernels in~$\mathbb R^n$ and later on generalized
by~B.~Fuglede~\cite{F1} for general kernels in a locally compact
space~$\mathbf X$. However, those two theories~--- for sets and, on
the other hand, condensers~--- are drastically different. To
illustrate this, it is enough to note that, in the noncompact case,
the main mi\-ni\-mum-pro\-blem of the theory of interior capacities
of condensers is in general {\it unsolvable\/}, and this phenomenon
occurs even under very natural assumptions (e.\,g., for the
Newtonian, Green, or Riesz kernels in~$\mathbb R^n$, $n\geqslant2$,
and closed condensers of finitely many plates); compare
with~\cite{Car, F1}. Necessary and sufficient conditions for the
problem to be solvable have been given in~\cite{Z4,Z6}; see~Sec.~5.1
below for a brief survey.

Therefore it was particularly interesting to find statements of
variational problems {\it dual\/} to the main minimum-problem of the
theory of interior capacities of condensers, but in contrast to the
last one, now {\it always solvable}~--- e.\,g., even for nonclosed,
unbounded condensers of infinitely many plates. (When speaking on
duality of variational problems, we mean their extremal values to be
equal.)

In all that follows, $\mathbf X$ denotes a locally compact Hausdorff
space, and $\mathfrak M=\mathfrak M(\mathbf X)$ the linear space of
all real-valued Radon measures~$\nu$ on~$\mathbf X$ equipped with
the {\it vague\/} ($={}${\it weak}$*$) topology, i.\,e., the
topology of pointwise convergence on the class $\mathbf C_0(\mathbf
X)$ of all real-valued continuous functions on~$\mathbf X$ with
compact support.

A {\it kernel\/} $\kappa$ on $\mathbf X$ is meant to be a lower
semicontinuous function $\kappa:\mathbf X\times\mathbf
X\to(-\infty,\infty]$. In order to avoid certain difficulties, we
follow~\cite{F1} in assuming that $\kappa\geqslant0$ unless the
space~$\mathbf X$ is compact.

The {\it energy\/} and the {\it potential\/} of a measure
$\nu\in\mathfrak M$ with respect to a kernel~$\kappa$ are defined by
$$\kappa(\nu,\nu):=\int\kappa(x,y)\,d(\nu\otimes\nu)(x,y)$$
and $$\kappa(x,\nu):=\int\kappa(x,y)\,d\nu(y),\quad x\in\mathbf X,$$
respectively, provided the corresponding integral above is well
defined (as a finite number or $\pm\infty$). Let $\mathcal E$ denote
the set of all $\nu\in\mathfrak M$ with
$-\infty<\kappa(\nu,\nu)<\infty$.

In the present study we shall be concerned with minimal energy
problems over certain subclasses of~$\mathcal E$, properly chosen.
For all positive definite kernels satisfying B.~Fuglede's condition
of consistency between the strong and the vague topologies
on~$\mathcal E$ (see~Sec.~2 below), those variational problems are
shown to be {\it dual\/} to the main minimum-pro\-blem of the theory
of interior capacities of con\-den\-sers (and hence providing some
new {\it equivalent\/} definitions of the capacity), but now {\it
always solvable\/}. See~Theorems~2\,--\,4 and Corollaries~10,~12.
Their solutions provide a natural generalization of the well-known
notion of interior capacitary distributions associated with a set
(see~\cite{F1}). We give a description of those solutions, establish
statements on their uniqueness and continuity, and point out their
characteristic properties; see~Sec.~7\,--\,10. The results obtained
hold true, e.\,g., for the Newtonian, Green or Riesz kernels
in~$\mathbb R^n$, $n\geqslant2$, as well as for the restriction of
the logarithmic kernel in~$\mathbb R^2$ to an open unit
ball.\bigskip

{\bf\large 2. Preliminaries: topologies, consistent and perfect
kernels}\nopagebreak\medskip

Recall that a measure $\nu\geqslant0$ is said to be {\it
concentrated\/} on~$E$, where $E$ is a subset of~$\mathbf X$, if the
complement $\complement E:=\mathbf X\setminus E$ is locally
$\nu$-negligible; or, equivalently, if $E$ is $\nu$-measurable and
$\nu=\nu_E$, where $\nu_E$ denotes the trace of~$\nu$ upon~$E$.

Let $\mathfrak M^+(E)$ be the convex cone of all nonnegative
measures concentrated on~$E$, and $\mathcal E^+(E):=\mathfrak
M^+(E)\cap\mathcal E$. We also write $\mathfrak M^+:=\mathfrak
M^+(\mathbf X)$ and $\mathcal E^+:=\mathcal E^+(\mathbf X)$.

From now on, the kernel under consideration is always assumed to be
{\it positive definite}, which means that it is symmetric (i.\,e.,
$\kappa(x,y)=\kappa(y,x)$ for all $x,\,y\in\mathbf X$) and the
energy $\kappa(\nu,\nu)$, $\nu\in\mathfrak M$, is nonnegative
whenever defined. Then $\mathcal E$ is known to be a pre-Hil\-bert
space with the scalar product
$$
\kappa(\nu_1,\nu_2):=\int\kappa(x,y)\,d(\nu_1\otimes\nu_2)(x,y)
$$
and the seminorm $\|\nu\|:=\sqrt{\kappa(\nu,\nu)}$; see~\cite{F1}. A
(positive definite) kernel is called {\it strictly positive
definite\/} if the seminorm $\|\cdot\|$ is a norm.

A measure $\nu\in\mathcal E$ is said to be {\it equivalent in\/}
$\mathcal E$ to a given $\nu_0\in\mathcal E$ if $\|\nu-\nu_0\|=0$;
the equivalence class, consisting of all those~$\nu$, will be
denoted by~$\left[\nu_0\right]_\mathcal E$.

In addition to the {\it strong\/} topology on~$\mathcal E$,
determined by the above seminorm~$\|\cdot\|$, it is often useful to
consider the {\it weak\/} topology on~$\mathcal E$, defined by means
of the seminorms $\nu\mapsto|\kappa(\nu,\mu)|$, $\mu\in\mathcal E$
(see~\cite{F1}). The Cauchy-Schwarz inequality
$$
|\kappa(\nu,\mu)|\leqslant\|\nu\|\,\|\mu\|,\quad\nu,\,\mu\in\mathcal
E,
$$
implies immediately that the strong topology on $\mathcal E$ is
finer than the weak one.

In \cite{F1}, B.~Fuglede introduced the following two properties of
{\it consistency\/} between the induced strong, weak, and vague
topologies on~$\mathcal E^+$:\smallskip

($C$) \ {\it Every strong Cauchy net in $\mathcal E^+$ converges
strongly to every its vague cluster point;}\smallskip

($CW$) \ {\it Every strongly bounded and vaguely convergent net in
$\mathcal E^+$ converges weakly to the vague limit;}\smallskip

\noindent in \cite{F2}, the properties ($C$) and ($CW$) were shown
to be {\it equivalent\/}.\medskip

{\bf Definition 1.} Following B.~Fuglede, we call a kernel~$\kappa$
{\it consistent} if it satisfies either of the properties~($C$)
and~($CW$), and {\it perfect\/} if, in addition, it is strictly
positive definite.\medskip

{\bf Remark 1.} One has to consider {\it nets\/} or {\it filters\/}
in~$\mathfrak M^+$ instead of sequences, for the vague topology in
general does not satisfy the first axiom of countability. We follow
Moore's and Smith's theory of convergence, based on the concept of
nets (see~\cite{MS}; cf.~also~\cite[Chap.~0]{E2} and
\cite[Chap.~2]{K}). However, if a locally compact space~$\mathbf X$
is metrizable and countable at infinity, then $\mathfrak M^+$
satisfies the first axiom of countability
(see~\cite[Lemma~1.2.1]{F1}) and the use of nets may be
avoided.\medskip

{\bf Theorem 1}~\cite{F1}. {\it A kernel $\kappa$ is perfect if and
only if $\mathcal E^+$ is strongly complete and the strong topology
on~$\mathcal E^+$ is finer than the vague one.}\medskip

{\bf Examples.} In $\mathbb R^n$, $n\geqslant 3$, the Newtonian
kernel $|x-y|^{2-n}$ is perfect~\cite{Car}. So are the Riesz kernel
$|x-y|^{\alpha-n}$, $0<\alpha<n$, in~$\mathbb R^n$, $n\geqslant2$
(see~\cite{D1, D2}), and the logarithmic kernel $-\log\,|x-y|$ in
$\mathbb R^2$, restricted to an open unit ball~\cite{L}.
Furthermore, if $D$ is an open set in~$\mathbb R^n$, $n\geqslant 2$,
and its generalized Green function~$g_D$ exists (see,
e.\,g.,~\cite[Th.~5.24]{HK}), then the Green kernel~$g_D$ is perfect
as well~\cite{E1}.\medskip

{\bf Remark 2.} As is seen from Theorem~1, the concept of consistent
or perfect kernels is an efficient tool in minimal energy problems
over classes of {\it nonnegative\/} measures with finite energy.
Indeed, the theory of capacities of {\it sets\/} has been developed
in~\cite{F1} exactly for those kernels. We shall show below that
this concept is still efficient in minimal energy problems over
classes of {\it signed\/} measures associated with a {\it
condenser}. This is guaranteed by a theorem on the strong
completeness of proper subspaces of~$\mathcal E$, to be stated in
Sec.~11 below.\bigskip

{\bf\large 3. Condensers of countably many plates. Measures
associated with a condenser; their energies and
potentials}\nopagebreak\medskip

{\bf 3.1.} Let $I^+$ and $I^-$ be countable (finite or infinite)
disjoint sets of indices~$i\in\mathbb N$, the latter being allowed
to be empty, and let $I$ denote their union. Assume that to every
$i\in I$ there corresponds a nonempty set~$A_i\subset\mathbf
X$.\medskip

{\bf Definition 2.} A collection $\mathcal A=(A_i)_{i\in I}$ is
called an $(I^+,I^-)$-{\it condenser\/} (or simply a {\it
condenser\/}) in~$\mathbf X$ if every compact subset of~$\mathbf X$
might have points in common with only a finite number of~$A_i$ and,
moreover,
\begin{equation}
\overline{A_i}\cap\overline{A_j}=\varnothing\quad\mbox{for all \ }
i\in I^+, \ j\in I^-. \label{non}
\end{equation}

The sets $A_i$, $i\in I^+$, and $A_j$, $j\in I^-$, are said to be
the {\it positive\/} and, respectively, the {\it negative plates\/}
of an $(I^+,I^-)$-condenser $\mathcal A=(A_i)_{i\in I}$. Note that
any two equal-sign\-ed plates of a condenser might intersect each
other (or even coincide).

Given $I^+$ and $I^-$, let $\mathfrak C=\mathfrak C(I^+,I^-)$ be the
class of all $(I^+,I^-)$-condensers in~$\mathbf X$. A condenser
$\mathcal A\in\mathfrak C$ is called {\it closed\/} or {\it
compact\/} if all $A_i$, $i\in I$, are closed or, respectively,
compact. Similarly, we call it {\it universally measurable\/} if all
the plates are universally measurable~--- that is, measurable with
respect to every $\nu\in\mathfrak M^+$. Next, $\mathcal
A=(A_i)_{i\in I}$ is said to be {\it finite\/} if so is~$I$.

Given $\mathcal A=(A_i)_{i\in I}$, write $\overline{\mathcal
A}:=(\,\overline{A_i}\,)_{i\in I}$. Then, due to~(\ref{non}),
$\overline{\mathcal A}$ is a (closed) $(I^+,I^-)$-condenser. In the
sequel, also the following notation will be required:
$$
A:=\bigcup_{i\in I}\,A_i,\qquad A^+:=\bigcup_{i\in I^+}\,A_i,\qquad
A^-:=\bigcup_{i\in I^-}\,A_i.
$$
Note that both $A^+$ and~$A^-$ might be noncompact even for a
compact~$\mathcal A$.\medskip

{\bf 3.2.} With the preceding notation, write
$$ \alpha_i:=\left\{
\begin{array}{rll} +1 & \mbox{if} & i\in I^+,\\ -1 & \mbox{if} & i\in
I^-.\\ \end{array} \right.
$$
Given $\mathcal A\in\mathfrak C$, let $\mathfrak M(\mathcal A)$
consist of all (finite or infinite) {\it linear combinations\/}
$$
\mu:=\sum_{i\in I}\,\alpha_i\mu^i,\quad\mbox{where \
}\mu^i\in\mathfrak M^+(A_i).
$$
Any two $\mu_1$ and $\mu_2$ in $\mathfrak M(\mathcal A)$,
$$\mu_1=\sum_{i\in
I}\,\alpha_i\mu_1^i\quad\mbox{and}\quad\mu_2=\sum_{i\in
I}\,\alpha_i\mu_2^i,
$$
are regarded to be {\it identical\/} ($\mu_1\equiv\mu_2$) if and
only if $\mu_1^i=\mu_2^i$ for all $i\in I$. Observe that, under the
relation of identity thus defined, the following correspondence
between $\mathfrak M(\mathcal A)$ and the Cartesian product
$\prod_{i\in I}\mathfrak M^+(A_i)$ is one-to-one:
$$\mathfrak
M(\mathcal A)\ni\mu\mapsto(\mu^i)_{i\in I}\in \prod_{i\in
I}\mathfrak M^+(A_i).$$
We call $\mu\in\mathfrak M(\mathcal A)$ a
{\it measure associated with}~$\mathcal A$, and $\mu^i$, $i\in I$,
its $i$-{\it coordinate}.

For measures associated with a condenser, it is therefore natural to
introduce the following concept of convergence, actually
corresponding to the vague convergence by coordinates. Let $S$
denote a directed set of indices, and let $\mu_s$, $s\in S$, and
$\mu_0$ be given elements of the class~$\mathfrak
M(\,\overline{\mathcal A}\,)$.\medskip

{\bf Definition 3.} A net $(\mu_s)_{s\in S}$ is said to converge to
$\mu_0$ $\mathcal A$-{\it va\-gue\-ly\/} if
$$
\mu^i_s\to\mu_0^i\quad\mbox{vaguely for all \ } i\in I.
$$

Then $\mathfrak M(\,\overline{\mathcal A}\,)$, equipped with the
topology of $\mathcal A$-vague convergence, and the product space
$\prod_{i\in I}\mathfrak M^+(\,\overline{A_i}\,)$ become
homeomorphic. Since the space $\mathfrak M(\mathbf X)$ is Hausdorff,
so are both $\mathfrak M(\,\overline{\mathcal A}\,)$ and
$\prod_{i\in I}\mathfrak M^+(\,\overline{A_i}\,)$ (see,
e.\,g.,~\cite[Chap.~3, Th.~5]{K}).

Similarly, a set $\mathfrak F\subset\mathfrak M(\,\overline{\mathcal
A}\,)$ is called $\mathcal A$-{\it vaguely bounded\/} if all its
$i$-projections are vaguely bounded~--- that is, if for all
$\varphi\in\mathbf C_0(\mathbf X)$ and $i\in I$,
$$\sup_{\mu\in\mathfrak F}\,|\mu^i(\varphi)|<\infty.$$

{\bf Lemma 1.} {\it If $\mathfrak F\subset\mathfrak
M(\,\overline{\mathcal A}\,)$ is bounded and closed in the $\mathcal
A$-vague topology, then it is $\mathcal A$-vaguely compact.}\medskip

{\bf Proof.} Since by~\cite[Chap.~III, \S~2, Prop.~9]{B2} any
vaguely bounded and closed part of~$\mathfrak M$ is vaguely compact,
the lemma follows immediately from Tychonoff's theorem on the
product of compact spaces (see, e.\,g.,~\cite[Chap.~5,
Th.~13]{K}).\hfill\phantom{B}\hfill$\Box$\medskip

{\bf 3.3.} Fix a linear combination $\mu\in\mathfrak M(\mathcal A)$.
Since each compact subset of~$\mathbf X$ might intersect with only
finite number of~$A_i$, $i\in I$, for every $\varphi\in\mathbf
C_0(\mathbf X)$ only finite number of~$\mu^i(\varphi)$, $i\in I$,
are nonzero. This yields that to every $\mu\in\mathfrak M(\mathcal
A)$ there corresponds a unique Radon measure~$R\mu$ such that
$$
R\mu(\varphi)=\sum_{i\in I}\,\alpha_i\mu^i(\varphi)\quad\mbox{for
all \ }\varphi\in\mathbf C_0(\mathbf X);
$$
its positive and negative parts in Jordan's decomposition, $R\mu^+$
and~$R\mu^-$, can be written in the form
$$
R\mu^+=\sum_{i\in I^+}\,\mu^i,\qquad R\mu^-=\sum_{i\in I^-}\,\mu^i.
$$

Of course, the mapping~$R:\mathfrak M(\mathcal A)\to\mathfrak M$
thus defined is in general non-injective, i.\,e., one may choose
$\mu'\in\mathfrak M(\mathcal A)$ so that $\mu'\not\equiv\mu$, while
$R\mu'=R\mu$. (It would be injective if all~$A_i$, $i\in I$, were
mutually disjoint.) We shall call $\mu,\,\mu'\in\mathfrak M(\mathcal
A)$ {\it equivalent in\/}~$\mathfrak M(\mathcal A)$, and write
$\mu\cong\mu'$, whenever their $R$-images coincide.\medskip

{\bf Lemma 2.} {\it The $\mathcal A$-vague convergence of
$(\mu_s)_{s\in S}$ to~$\mu_0$ implies the vague convergence of
$(R\mu_s)_{s\in S}$ to~$R\mu_0$.}\medskip

{\bf Proof.} This is obvious in view of the fact that the support of
any $\varphi\in\mathbf C_0(\mathbf X)$ may have points in common
with only a finite number of $\overline{A_i}$, $i\in
I$.\hfill\phantom{B}\hfill$\Box$\medskip

{\bf Remark 3.} The statement of Lemma~2 in general can not be
inverted. However, if all $\overline{A_i}$, $i\in I$, are mutually
disjoint, then the vague convergence of $(R\mu_s)_{s\in S}$
to~$R\mu_0$ implies the $\mathcal A$-vague convergence of
$(\mu_s)_{s\in S}$ to~$\mu_0$. This can be seen by using the
Tietze-Urysohn extension theorem (see,
e.\,g.,~\cite[Th.~0.2.13]{E2}).\medskip

{\bf 3.4.} We next proceed to define energies and potentials of
$\mu\in\mathfrak M(\mathcal A)$. A proper definition is based on the
mapping $R:\mathfrak M(\mathcal A)\to\mathfrak M$ and the following
assertion.\medskip

{\bf Lemma 3.} {\it Fix $\mu\in\mathfrak M(\mathcal A)$ and a lower
semicontinuous function~$\psi$ on~$\mathbf X$ such that
$\psi\geqslant0$ unless $\mathbf X$~is compact. If the integral
$\int\psi\,d R\mu$ is well defined, then
\begin{equation}
\int\psi\,dR\mu=\sum_{i\in
I}\,\alpha_i\int\psi\,d\mu^i,\label{lemma11}
\end{equation}
and it is finite if and only if the series on the right is
absolutely convergent.}\medskip

{\bf Proof.} We can certainly assume $\psi$ to be nonnegative, for
if not, we replace $\psi$ by a function~$\psi'$ obtained by adding
to~$\psi$ a suitable constant~$c>0$:
$$
\psi'(x):=\psi(x)+c\geqslant0,
$$
which is always possible since a lower semicontinuous function is
bounded from below on a compact space. Hence, for every $N\in\mathbb
N$,
$$\int\psi\,d R\mu^+\geqslant\sum_{i\in I^+, \ i\leqslant N}\,\int\psi\,d\mu^i.$$
On the other hand, the sum of $\mu^i$ over all $i\in I^+$ that do
not exceed~$N$ approaches~$R\mu^+$ vaguely as $N\to\infty$;
consequently (see, e.\,g.,~\cite{F1})
$$\int\psi\,d R\mu^+\leqslant\lim_{N\to\infty}\,\sum_{i\in I^+, \ i\leqslant N}\,\int\psi\,d\mu^i.$$
Combining the last two inequalities and then letting $N$ tend
to~$\infty$ yields
$$\int\psi\,d R\mu^+=\sum_{i\in
I^+}\,\int\psi\,d\mu^i.$$  Since the same holds true for $R\mu^-$
and~$I^-$ instead of~$R\mu^+$ and~$I^+$, respectively, the lemma
follows.\hfill\phantom{B}\hfill$\Box$\medskip

{\bf Corollary 1.} {\it Given $\mu,\,\mu_1\in\mathfrak M(\mathcal
A)$ and $x\in\mathbf X$, then
\begin{align}
\kappa(x,R\mu)&=\sum_{i\in
I}\,\alpha_i\kappa(x,\mu^i),\label{poten}\\[4pt]
\kappa(R\mu,R\mu_1)&=\sum_{i,j\in
I}\,\alpha_i\alpha_j\kappa(\mu^i,\mu_1^j),\label{mutual}
\end{align}
each of the identities being understood in the sense that its
right-hand side is well defined whenever so is the left-hand one and
then they coincide. Furthermore, the left-hand side
in\/}~(\ref{poten}) {\it or in\/}~(\ref{mutual}) {\it is finite if
and only if the corresponding series on the right absolutely
converges.}\medskip

{\bf Proof.} Relation (\ref{poten}) is a direct consequence
of~(\ref{lemma11}), while (\ref{mutual}) follows from Fubini's
theorem and Lemma~3 on account of the fact that $\kappa(x,\nu)$,
where $\nu\in\mathfrak M^+$ is given, is lower semicontinuous
on~$\mathbf X$ (see,
e.\,g.,~\cite{F1}).\hfill\phantom{B}\hfill$\Box$\medskip

{\bf Definition 4.} Given $\mu,\,\mu_1\in\mathfrak M(\mathcal A)$,
then
$$\kappa(x,\mu):=\kappa(x,R\mu)$$
is called the value of the {\it potential\/} of~$\mu$ at a point
$x\in\mathbf X$, and
$$\kappa(\mu,\mu_1):=\kappa(R\mu,R\mu_1)$$ the {\it mutual
energy\/} of~$\mu$ and~$\mu_1$~--- of course, provided the
right-hand side of the corresponding relation is well defined. For
$\mu\equiv\mu_1$ we get the {\it energy\/} $\kappa(\mu,\mu)$
of~$\mu$; i.\,e., if $\kappa(R\mu,R\mu)$ is well defined, then
\begin{equation}
\kappa(\mu,\mu):=\kappa(R\mu,R\mu)=\sum_{i,j\in
I}\,\alpha_i\alpha_j\kappa(\mu^i,\mu^j).\label{ener}\end{equation}

{\bf Corollary~2.} {\it For $\mu\in\mathfrak M(\mathcal A)$ to be of
finite energy, it is necessary and sufficient that so be all
$\mu^i$, $i\in I$, and}
$$\sum_{i\in I}\,\|\mu^i\|^2<\infty.$$

{\bf Proof.} This follows immediately from (\ref{ener}) and
Corollary~1 due to the inequality
$$2\kappa(\nu_1,\nu_2)\leqslant\|\nu_1\|^2+\|\nu_2\|^2,\quad \nu_1,\,\nu_2\in\mathcal E\vspace{-8mm}.$$
\hfill$\Box$\medskip

{\bf Remark 4.} Given $\mu\in\mathfrak M(\mathcal A)$, then the
series in~(\ref{ener}) actually defines the energy of the vector
measure $(\mu^i)_{i\in I}$ relative to the infinite interaction
matrix of the form $(\alpha_i\alpha_j)_{i,j\in I}$\,; compare
with~\cite{GR}, \cite[Chap.~5, \S~4]{NS}. Our approach, however, is
based on the fact that, due to the specific interaction matrix, this
value can also be obtained as the energy of the corresponding Radon
measure~$R\mu$.\medskip

{\bf Remark 5.} Since we make no difference between $\mu\in\mathfrak
M(\mathcal A)$ and $R\mu$ when dealing with their energies or
potentials, we shall sometimes call a measure associated with a
condenser simply a {\it measure\/}~--- certainly, if this causes no
confusion.\medskip

{\bf 3.5.} Let $\mathcal E(\mathcal A)$ consist of all
$\mu\in\mathfrak M(\mathcal A)$ of finite energy
$\kappa(\mu,\mu)=:\|\mu\|^2$. Since $\mathfrak M(\mathcal A)$ forms
a convex cone, it is seen from Corollary~2 that so does $\mathcal
E(\mathcal A)$. Let us treat $\mathcal E(\mathcal A)$ as a {\it
semimetric space\/} with the semimetric
\begin{equation}
\|\mu_1-\mu_2\|:=\|R\mu_1-R\mu_2\|,\quad\mu_1,\,\mu_2\in\mathcal
E(\mathcal A);\label{seminorm}
\end{equation}
then $\mathcal E(\mathcal A)$ and its $R$-image become isometric.
The topology on $\mathcal E(\mathcal A)$ defined by means of the
semimetric~(\ref{seminorm}) will be called {\it strong\/}. Two
elements of $\mathcal E(\mathcal A)$, $\mu_1$ and~$\mu_2$, are
called {\it equivalent in\/}~$\mathcal E(\mathcal A)$ if
$\|\mu_1-\mu_2\|=0$. If, in addition, the kernel~$\kappa$ is assumed
to be strictly positive definite, then the equivalence in~$\mathcal
E(\mathcal A)$ implies that in~$\mathfrak M(\mathcal A)$, namely
then $\mu_1\cong\mu_2$.\bigskip

{\bf\large 4. Interior capacities of condensers; elementary
properties}\nopagebreak\medskip

{\bf 4.1.} Let $\mathcal H$ be a set in the pre-Hilbert
space~$\mathcal E$ or in the semimetric space~$\mathcal E(\mathcal
A)$, an $(I^+,I^-)$-con\-den\-ser~$\mathcal A$ being given. In
either case, let us introduce the quantity
$$
\|\mathcal H\|^2:=\inf_{\nu\in\mathcal H}\,\|\nu\|^2,
$$
interpreted as $+\infty$ if $\mathcal H$ is empty. If $\|\mathcal
H\|^2<\infty$, one can consider the variational problem on the
existence of $\lambda=\lambda(\mathcal H)\in\mathcal H$ with minimal
energy
$$
\|\lambda\|^2=\|\mathcal H\|^2;
$$
such a problem will be referred to as the $\mathcal H$-{\it
problem\/}. The $\mathcal H$-problem is called {\it solvable\/} if a
minimizing measure $\lambda(\mathcal H)$ exists.

The following elementary lemma is a slight generalization of
\cite[Lemma~4.1.1]{F1}.\medskip

{\bf Lemma 4.} {\it Suppose $\mathcal H$ is convex and
$\lambda=\lambda(\mathcal H)$ exists. Then for any $\nu\in\mathcal
H$,}
\begin{equation}
\|\nu-\lambda\|^2\leqslant\|\nu\|^2-\|\lambda\|^2.\label{lemma1}
\end{equation}

{\bf Proof.} Assume $\mathcal H\subset\mathcal E$. For every
$t\in[0,\,1]$, the measure $\mu:=(1-t)\lambda+t\nu$ belongs
to~$\mathcal H$, and therefore $\|\mu\|^2\geqslant\|\lambda\|^2$.
Evaluating~$\|\mu\|^2$ and then letting $t$ tend to zero, we get
$\kappa(\nu,\lambda)\geqslant\|\lambda\|^2$, and (\ref{lemma1})
follows (see~\cite{F1}).

Suppose now $\mathcal H\subset\mathcal E(\mathcal A)$. Then
$R\mathcal H:=\{R\nu:\ \nu\in\mathcal H\}$ is a convex subset
of~$\mathcal E$, while $R\lambda$ is a minimizer in the $R\mathcal
H$-problem. What has just been shown therefore yields
$$\|R\nu-R\lambda\|^2\leqslant\|R\nu\|^2-\|R\lambda\|^2,$$
which gives (\ref{lemma1}) when combined
with~(\ref{seminorm}).\hfill\phantom{B}\hfill$\Box$\medskip

We shall be concerned with the $\mathcal H$-problem for various
specific~$\mathcal H$ related to the notion of {\it interior
capacity\/} of an $(I^+,I^-)$-con\-den\-ser (in particular, of a
set); see~Sec.~4.2 and~Sec.~7 below for the definitions.\medskip

{\bf 4.2.} Fix a continuous function $g:\mathbf X\to(0,\infty)$ and
a numerical vector $a=(a_i)_{i\in I}$ with $a_i>0$, $i\in I$. Given
an $(I^+,I^-)$-condenser~$\mathcal A$ in~$\mathbf X$, write
$$
\mathfrak M^+(A_i,a_i,g):=\Bigl\{\nu\in\mathfrak M^+(A_i): \ \int
g\,d\nu=a_i\Bigr\},
$$
and let $\mathfrak M(\mathcal A,a,g)$ consist of all
$\mu\in\mathfrak M(\mathcal A)$ with $\mu^i\in\mathfrak
M^+(A_i,a_i,g)$ for all $i\in I$. Given a kernel~$\kappa$, also
write
$$\mathcal E^+(A_i,a_i,g):=\mathfrak M^+(A_i,a_i,g)\cap\mathcal
E,\quad\mathcal E(\mathcal A,a,g):=\mathfrak M(\mathcal
A,a,g)\cap\mathcal E(\mathcal A).$$

{\bf Definition 5.} We shall call the value
\begin{equation}
{\rm cap}\,\mathcal A:={\rm cap}\,(\mathcal
A,a,g):=\frac{1}{\|\mathcal E(\mathcal A,a,g)\|^2} \label{def}
\end{equation}
the ({\it interior}) {\it capacity\/} of an $(I^+,I^-)$-condenser
$\mathcal A$ (with respect to~$\kappa$, $a$, and~$g$).\medskip

Here and in the sequel, we adopt the convention that $1/0=+\infty$.
It follows immediately from the positive definiteness of the kernel
that
$$
0\leqslant{\rm cap}\,(\mathcal A,a,g)\leqslant\infty.
$$

{\bf Remark 6.} If $I$ is a singleton, then any
$(I^+,I^-)$-con\-den\-ser consists of just one set, say~$A_1$. If
moreover $g=1$ and $a_1=1$, then the notion of interior capacity of
a condenser, defined above, certainly reduces to the notion of
interior capacity of a set (see~\cite{F1}). We denote it
by~$C(\,\cdot\,)$ as well, i.\,e., $C(A_1):=\|\mathcal
E^+(A_1,1,1)\|^{-2}$.\medskip

{\bf Remark 7.} In the case of the Newtonian kernel $|x-y|^{-1}$ in
$\mathbb R^3$, the notion of capacity of a condenser~$\mathcal A$
has an evident electrostatic interpretation. In the framework of the
corresponding electrostatics problem, the function~$g$ serves as a
characteristic of nonhomogeneity of the conductors $A_i$, $i\in
I$.\medskip

{\bf 4.3.} On $\mathfrak C=\mathfrak C(I^+,I^-)$, it is natural to
introduce an ordering relation~$\prec$ by declaring $\mathcal
A'\prec\mathcal A$ to mean that $A_i'\subset A_i$ for all $i\in I$.
Here, $\mathcal A'=(A_i')_{i\in I}$. Then ${\rm
cap}\,(\,\cdot\,,a,g)$ is a nondecreasing function of a condenser,
namely
\begin{equation}
{\rm cap}\,(\mathcal A',a,g)\leqslant{\rm cap}\,(\mathcal
A,a,g)\quad\mbox{whenever \ }\mathcal A'\prec\mathcal A.
\label{increas'}
\end{equation}

Given $\mathcal A\in\mathfrak C$, denote by $\{\mathcal K\}_\mathcal
A$ the increasing ordered family of all compact condensers $\mathcal
K=(K_i)_{i\in I}\in\mathfrak C$ such that $\mathcal K\prec\mathcal
A$.\medskip

{\bf Lemma 5.} {\it If $\mathcal K$ ranges over $\{\mathcal
K\}_\mathcal A$, then}
\begin{equation}
{\rm cap}\,(\mathcal A,a,g)=\lim_{\mathcal K\uparrow\mathcal A}\,
{\rm cap}\,(\mathcal K,a,g).\label{cont}
\end{equation}

{\bf Proof.} We can certainly assume ${\rm cap}\,(\mathcal A,a,g)$
to be nonzero, since otherwise (\ref{cont})~follows at once
from~(\ref{increas'}). Then the set $\mathcal E(\mathcal A,a,g)$
must be nonempty; fix~$\mu$, one of its elements. Given $\mathcal
K\in\{\mathcal K\}_\mathcal A$ and $i\in I$, let $\mu^i_\mathcal K$
denote the trace of~$\mu^i$ upon~$K_i$, i.\,e., $\mu^i_\mathcal
K:=\mu_{K_i}^i$. Applying Lemma~1.2.2 from~\cite{F1}, we conclude
that
\begin{align}
\int g\,d\mu^i&=\lim_{\mathcal K\uparrow\mathcal A}\,\int
g\,d\mu_{\mathcal K}^i,\hspace{-3cm} &i&\in I,\label{w}\\
\kappa(\mu^i,\mu^j)&=\lim_{\mathcal K\uparrow\mathcal
A}\,\kappa(\mu_{\mathcal K}^i,\mu_{\mathcal K}^j),\hspace{-3cm}
&i,\,j&\in I.\label{ww}
\end{align}

Fix $\varepsilon>0$. It follows from~(\ref{w}) and~(\ref{ww}) that
for every $i\in I$ one can choose a compact set $K_i^0\subset A_i$
so that
\begin{equation}
\frac{a_i}{\int
g\,d\mu^i_{K_i^0}}<1+\varepsilon\,i^{-2},\label{unif2}
\end{equation}
\begin{equation}
\bigl|\|\mu^i\|^2-\|\mu^i_{K_i^0}\|^2\bigr|<\varepsilon^2i^{-4}.\label{unif1}
\end{equation}
Having denoted $\mathcal K^0:=(K_i^0)_{i\in I}$, for every $\mathcal
K\in\{\mathcal K\}_\mathcal A$ that follows~$\mathcal K^0$ we
therefore have $\int g\,d\mu_\mathcal K^i\ne0$ and
$$
\hat{\mu}_\mathcal K:=\sum_{i\in I}\,\frac{\alpha_ia_i}{\int
g\,d\mu_\mathcal K^i}\,\mu_\mathcal K^i\in\mathcal E(\mathcal
K,a,g),
$$
the finiteness of the energy being obtained from~(\ref{unif1})
and~Corollary~2. This yields
\begin{equation}
\|\hat{\mu}_\mathcal K\|^2\geqslant\|\mathcal E(\mathcal
K,a,g)\|^2.\label{www}\end{equation}

We next proceed to show that
\begin{equation}
\|\mu\|^2=\lim_{\mathcal K\uparrow\mathcal A}\,\|\hat{\mu}_\mathcal
K\|^2.\label{4w}\end{equation} To this end, it can be assumed that
$\kappa\geqslant0$; for if not, then $\mathcal A$ must be finite
since $\mathbf X$ is compact, and (\ref{4w}) follows from~(\ref{w})
and~(\ref{ww}) when substituted into~(\ref{ener}). Therefore, for
every~$\mathcal K$ that follows~$\mathcal K_0$ and every $i\in I$ we
get
\begin{equation}
\|\mu^i_\mathcal K\|\leqslant\|\mu^i\|\leqslant\|R\mu^++R\mu^-\|,
\label{unif3}
\end{equation}
\begin{equation}
\|\mu^i-\mu^i_{\mathcal K}\|<\varepsilon\,i^{-2},\label{uniff}
\end{equation}
the latter being clear from (\ref{unif1}) because of
$\kappa(\mu^i_{\mathcal K},\mu^i-\mu^i_{\mathcal K})\geqslant0$.
Also observe that, by~(\ref{ener}),
\begin{equation*}
\begin{split}
\bigl|&\|\mu\|^2-\|\hat{\mu}_\mathcal
K\|^2\bigr|\leqslant\sum_{i,j\in
I}\,\biggl|\kappa(\mu^i,\mu^j)-\frac{a_i}{\int g\,d\mu^i_\mathcal
K}\frac{a_j}{\int g\,d\mu^j_\mathcal K}\,\kappa(\mu_\mathcal
K^i,\mu_\mathcal K^j)\biggr|\\[7pt]
&{}\leqslant\sum_{i,j\in I}\,\biggl[\kappa(\mu^i-\mu^i_\mathcal
K,\mu^j)+\kappa(\mu^i_\mathcal K,\mu^j-\mu^j_\mathcal
K)+\biggl(\frac{a_i}{\int g\,d\mu^i_\mathcal K}\frac{a_j}{\int
g\,d\mu^j_\mathcal K}-1\biggr)\kappa(\mu^i_\mathcal K,\mu^j_\mathcal
K)\biggr].
\end{split}
\end{equation*}
When combined with (\ref{unif2}), (\ref{unif3}), and~(\ref{uniff}),
this yields
$$\bigl|\|\mu\|^2-\|\hat{\mu}_\mathcal
K\|^2\bigr|\leqslant M\varepsilon\quad\mbox{for all \ }\mathcal
K\succ\mathcal K_0,$$ where $M$ is  finite and independent
of~$\mathcal K$, and the required relation~(\ref{4w}) follows.

Substituting (\ref{www}) into~(\ref{4w}), in view of the arbitrary
choice of $\mu\in\mathcal E(\mathcal A,a,g)$ we get
$$
\|\mathcal E(\mathcal A,a,g)\|^2\geqslant\lim_{\mathcal
K\uparrow\mathcal A}\,\|\mathcal E(\mathcal K,a,g)\|^2.
$$
Since the converse inequality is obvious from~(\ref{increas'}), the
proof is complete.\hfill\phantom{B}\hfill$\Box$\medskip

Let $\mathcal E^0(\mathcal A,a,g)$ denote the class of all
$\mu\in\mathcal E(\mathcal A,a,g)$ such that, for every $i\in I$,
the support~$S(\mu^i)$ of~$\mu^i$ is compact and contained
in~$A_i$.\medskip

{\bf Corollary 3.} {\it The capacity ${\rm cap}\,(\mathcal A,a,g)$
remains unchanged if the class $\mathcal E(\mathcal A,a,g)$ in its
definition is replaced by $\mathcal E^0(\mathcal A,a,g)$.  In other
words,}
$$
\|\mathcal E(\mathcal A,a,g)\|^2=\|\mathcal E^0(\mathcal A,a,g)\|^2.
$$

{\bf Proof.} We can certainly assume $\|\mathcal E(\mathcal
A,a,g)\|^2$ to be finite, since otherwise the corollary follows
immediately from $\mathcal E^0(\mathcal A,a,g)\subset\mathcal
E(\mathcal A,a,g)$. Then, by~(\ref{increas'}) and~(\ref{cont}), for
every $\varepsilon>0$ there exists a compact condenser $\mathcal
K\prec\mathcal A$ such that
$$\|\mathcal E(\mathcal K,a,g)\|^2\leqslant\|\mathcal E(\mathcal A,a,g)\|^2+\varepsilon.$$
This leads to the claimed assertion when combined with the relation
$$\|\mathcal E(\mathcal K,a,g)\|^2\geqslant\|\mathcal
E^0(\mathcal A,a,g)\|^2\geqslant\|\mathcal E(\mathcal
A,a,g)\|^2\vspace{-8mm}.$$\hfill$\Box$\medskip

{\bf 4.4.} Unless explicitly stated otherwise, in all that follows
it is assumed that
\begin{equation}
{\rm cap}\,(\mathcal A,a,g)>0;\label{nonzero1}
\end{equation}
see below for necessary and sufficient conditions for this to
occur.\medskip

{\bf Lemma 6.} {\it For\/} (\ref{nonzero1}) {\it to hold, it is
necessary and sufficient that either of the following three
equivalent conditions be satisfied\/}:\\[5pt]
(i)\quad$\mathcal E(\mathcal A,a,g)$ {\it is nonempty\/};\\[5pt]
(ii)\quad{\it there exist $\nu_i\in\mathcal E^+(A_i,a_i,g)$,
$i\in I$, such that\/} $\sum_{i\in I}\,\|\nu_i\|^2<\infty$;\\[5pt]
(iii)\quad$\sum_{i\in I}\,\|\mathcal
E^+(A_i,a_i,g)\|^2<\infty$.\medskip

{\bf Proof.} Indeed, the equivalency of (\ref{nonzero1}) and (i) is
obvious, while that of~(i) and~(ii) can be obtained directly from
Corollary~2. If (iii)~holds, then one can choose $\nu_i\in\mathcal
E^+(A_i,a_i,g)$, $i\in I$, so that $\|\nu_i\|^2<\|\mathcal
E^+(A_i,a_i,g)\|^2+i^{-2}$, and (ii) follows. Since (ii)~obviously
results~(iii), the proof is
complete.\hfill\phantom{B}\hfill$\Box$\medskip

{\bf Corollary 4.} {\it For\/} (\ref{nonzero1}) {\it to be
satisfied, it is necessary that
\begin{equation}
C(A_i)>0\quad\mbox{for all \ } i\in I. \label{cap2}
\end{equation}
If $\mathcal A$ is finite, then\/} (\ref{nonzero1}) {\it and\/}
(\ref{cap2}) {\it are actually equivalent.}\medskip

{\bf Proof.} For Lemma 6, (ii) to hold, it is necessary that, for
every $i\in I$, there exists a nonzero nonnegative measure of finite
energy concentrated on~$A_i$, which in turn is equivalent
to~(\ref{cap2}) by~\cite[Lemma~2.3.1]{F1}. Since the former
implication can obviously be inverted whenever $\mathcal A$ is
finite, the proof is complete.\hfill\phantom{B}\hfill$\Box$\medskip

Let $g_{\inf}$ and $g_{\sup}$ denote respectively the infimum and
the supremum of~$g$ over~$A$.\medskip

{\bf Corollary 5.} {\it Assume $0<g_{\inf}\leqslant
g_{\sup}<\infty$. Then\/} (\ref{nonzero1}) {\it holds if and only
if\/}
$$\sum_{i\in I}\,\frac{a_i^2}{C(A_i)}<\infty.$$

{\bf Proof.} Lemma~6, (iii) implies the corollary when combined with
the inequalities
\begin{equation}
\frac{a^2_i}{g_{\sup}^2\,C(A_i)}\leqslant\|\mathcal
E^+(A_i,a_i,g)\|^2\leqslant\frac{a^2_i}{g_{\inf}^2\,C(A_i)}\,,\quad
i\in I,\label{5w}
\end{equation}
to be proved below by reasons of homogeneity.

To verify (\ref{5w}), fix $i\in I$. One can certainly assume
$C(A_i)$ to be nonzero, for otherwise Corollary~4 with $I=\{i\}$
shows that each of the three parts in~(\ref{5w}) equals~$+\infty$.
Therefore, there exists $\theta\in\mathcal E^+(A_i,1,1)$. Since
$$\theta':=\frac{a_i\theta}{\int g\,d\theta}\in\mathcal
E^+(A_i,a_i,g),$$ we get
$$a_i^2\,\|\theta\|^2\geqslant g_{\inf}^2\,\|\theta'\|^2\geqslant
g_{\inf}^2\,\|\mathcal E^+(A_i,a_i,g)\|^2,$$ and the right-hand side
of~(\ref{5w}) is obtained by letting $\theta$ range over $\mathcal
E^+(A_i,1,1)$.

To verify the left-hand side, fix $\omega\in\mathcal
E^+(A_i,a_i,g)$. Then
$$0<a_i\,g_{\sup}^{-1}\leqslant\omega(\mathbf X)\leqslant
a_i\,g_{\inf}^{-1}<\infty.$$ Hence, $\omega(\mathbf
X)^{-1}\omega\in\mathcal E^+(A_i,1,1)$ and
$$\|\omega\|^2\geqslant\frac{a_i^2}{g_{\sup}^2}\,\|\mathcal E^+(A_i,1,1)\|^2.$$
In view of the arbitrary choice of $\omega\in\mathcal
E^+(A_i,a_i,g)$, this completes the
proof.\hfill\phantom{B}\hfill$\Box$\medskip

{\bf 4.5.} In the following assertion, providing necessary and
sufficient conditions for ${\rm cap}\,(\mathcal A,a,g)$ to be
finite, we assume that $g_{\inf}>0$.\medskip

{\bf Lemma 7.} {\it For ${\rm cap}\,(\mathcal A,a,g)$ to be finite,
it is necessary that
\begin{equation}
C(A_j)<\infty\quad\mbox{for some \ } j\in I.\label{j}
\end{equation}
This condition is also sufficient if it is additionally assumed that
$\sum_{i\in I}\,a_i<\infty$, $g_{\sup}<\infty$, $\mathcal A$ is
closed, while $\kappa$ is bounded from above on~$A^+\times A^-$ and
perfect.}\medskip

{\bf Proof.} Let ${\rm cap}\,\mathcal A<\infty$ and assume, on the
contrary, that
\begin{equation}
C(A_i)=\infty\quad\mbox{for all \ } i\in I.\label{lemma6}
\end{equation}
Given $\varepsilon>0$, then for every~$i$ one can choose
$\nu_i\in\mathcal E^+(A_i,1,1)$ with compact support so that
$$
\|\nu_i\|\leqslant\varepsilon a_i^{-1}i^{-2}\,g_{\inf}.
$$
Since then
$$
\nu:=\sum_{i\in I}\,\frac{\alpha_i a_i\nu_i}{\int
g\,d\nu_i}\in\mathcal E(\mathcal A,a,g)
$$
and
$$
\|\nu\|\leqslant\varepsilon\sum_{i\in I}\,i^{-2},
$$
we arrive at a contradiction by letting $\varepsilon$ tend to $0$.

Assume now all the conditions of the remaining part of the lemma to
be satisfied, and let (\ref{j}) be true~--- say, for $j\in I^+$.
Consider the (finite) condenser~$\mathcal B$ with the positive
plates~$B_1$ and~$B_2$ and the negative one~$B_3$, where
$$B_1:=A_j,\quad B_2:=\bigcup_{i\in I^+\setminus\{j\}}\,A_i,\quad
B_3:=A^-.$$ Also write $b:=(b_1,b_2,b_3)$, where
$$b_1:=a_j,\quad b_2:=\sum_{i\in I^+\setminus\{j\}}\,a_i,\quad b_3:=\sum_{i\in
I^-}\,a_i.$$ (If either of the sets $B_2$ and $B_3$ is empty, it
should just be dropped, as well as the corresponding coordinate
of~$b$.) Then for every $\mu\in\mathcal E(\mathcal A,a,g)$ there
exists $\nu\in\mathcal E(\mathcal B,b,g)$ such that $R\mu=R\nu$, and
therefore
$$
\|\mathcal E(\mathcal B,b,g)\|^2\leqslant\|\mathcal
E(\mathcal A,a,g)\|^2.
$$
Furthemore, Lemma~13 from~\cite{Z5} shows
that, under the stated assumptions, there exists $\zeta\in\mathcal
E(\mathcal B)$ such that $\int g\,d\zeta^1=b_1$ (hence,
$\zeta\not\equiv0$) and
$$\|\zeta\|^2=\|\mathcal E(\mathcal B,b,g)\|^2.$$
Since $\kappa$ is strictly positive definite, this implies that
$\|\mathcal E(\mathcal B,b,g)\|^2$ is nonzero. Hence, so is
$\|\mathcal E(\mathcal A,a,g)\|^2$, as was to be
proved.\hfill\phantom{B}\hfill$\Box$\bigskip

{\bf\large 5. On the solvability of the main
minimum-pro\-blem}\nopagebreak\medskip

{\bf 5.1.} Because of~(\ref{nonzero1}), we are naturally led to the
$\mathcal E(\mathcal A,a,g)$-{\it problem\/} (cf.~Sec.~4.1), i.\,e.,
the problem on the existence of $\lambda\in\mathcal E(\mathcal
A,a,g)$ with minimal energy
$$
\|\lambda\|^2=\|\mathcal E(\mathcal A,a,g)\|^2;
$$
the $\mathcal E(\mathcal A,a,g)$-problem might certainly be regarded
as the main minimum-problem of the theory of interior capacities of
condensers. The collection (possibly empty) of all minimizing
measures~$\lambda$ in this problem will be denoted by $\mathcal
S(\mathcal A,a,g)$.

If moreover ${\rm cap}\,(\mathcal A,a,g)<\infty$, let us look, as
well, at the $\mathcal E(\mathcal A,a\,{\rm cap}\,\mathcal
A,g)$-{\it prob\-lem\/}. By reasons of homogeneity, both the
$\mathcal E(\mathcal A,a,g)$- and the $\mathcal E(\mathcal A,a\,{\rm
cap}\,\mathcal A,g)$-problems are simultaneously either solvable or
unsolvable, and their extremal values are related to each other by
the following law:
\begin{equation}
\frac{1}{\|\mathcal E(\mathcal A,a,g)\|^2}=\|\mathcal E(\mathcal
A,a\,{\rm cap}\,\mathcal A,g)\|^2. \label{iden}
\end{equation}

Assume for a moment that $\mathcal A$ is compact. Since the mapping
$$\nu\mapsto\int g\,d\nu,\quad \nu\in\mathfrak M^+(K),$$
where $K\subset\mathbf X$ is a compact set, is vaguely continuous,
$\mathfrak M(\mathcal A,a,g)$ is compact in the $\mathcal A$-vague
topology. Therefore, if $\mathcal A$ is additionally assumed to be
finite, while $\kappa$ is continuous on~$A^+\times A^-$ (which, due
to~(\ref{non}), is always the case for either of the classical
kernels), then $\|\mu\|^2$ is $\mathcal A$-vaguely lower
semicontinuous on~$\mathcal E(\mathcal A)$, and the solvability of
both the problems immediately follows (cf.~\cite[Th.~2.30]{O}).

But if $\mathcal A$ is {\it noncompact\/}, then the class $\mathfrak
M(\mathcal A,a,g)$ is no longer $\mathcal A$-vaguely compact and the
problems become quite nontrivial. Moreover, it has recently been
shown by the author that, in the noncompact case, the problems are
in general {\it unsolvable\/} and this occurs even under very
natural assumptions (e.\,g., for the Newtonian, Green, or Riesz
kernels in~$\mathbb R^n$, $n\geqslant 2$, and finite, closed
condensers).

In particular, it was proved in~\cite{Z4} that, if $\mathcal A$ is
finite and closed, $\kappa$~is perfect, and bounded and continuous
on~$A^+\times A^-$, and satisfies the generalized maximum principle
(see, e.\,g.,~\cite[Chap.~VI]{L}), while $0<g_{\inf}\leqslant
g_{\sup}<\infty$, then for either of the $\mathcal E(\mathcal
A,a,g)$- and the $\mathcal E(\mathcal A,a\,{\rm cap}\,\mathcal
A,g)$-prob\-lems to be solvable for any vector~$a$, it is necessary
and sufficient that
$$
C(A_i)<\infty\quad\mbox{for all \ } i\in I.
$$
If moreover there exists $i_0\in I$ such that
$$
C(A_{i_0})=\infty,
$$
then both the problems are unsolvable for all $a=(a_i)_{i\in I}$
with $a_{i_0}$ large enough.

In~\cite[Th.~1]{Z6}, the last statement was sharpened. It was shown
that if, in addition to all the preceding assumptions, for all
$i\neq i_0$,
$$
C(A_i)<\infty\quad\mbox{and}\quad A_i\cap A_{i_0}=\varnothing,
$$
while $\kappa(\cdot,y)\to0$ (as $y\to\infty$) uniformly on compact
sets, then there exists a number $\Lambda_{i_0}\in[0,\infty)$ such
that the problems are unsolvable if and only if
$$a_{i_0}>\Lambda_{i_0}.$$

{\bf Remark 8.} It was actually shown in \cite{Z6} that
$$
\Lambda_{i_0}=\int g\,d\tilde{\lambda}^{i_0},
$$
where $\tilde{\lambda}$ is a minimizer (it exists) in the auxiliary
$\mathcal H$-problem for
$$
\mathcal H:=\left\{\mu\in\mathcal E(\mathcal
A):\quad\mu^i\in\mathcal E^+(A_i,a_i,g)\quad\mbox{for all \ } i\neq
i_0\right\}.
$$

{\bf Remark 9.} The mentioned results were actually obtained
in~\cite{Z4,Z6} for the energy evaluated in the presence of an
external field.\medskip

{\bf 5.2.} In view of the results reviewed in Sec.~5.1, it was
particularly interesting to find statements of variational problems
{\it dual\/} to the $\mathcal E(\mathcal A,a\,{\rm cap}\,\mathcal
A,g)$-problem (and hence providing new {\it equivalent\/}
definitions of ${\rm cap}\,\mathcal A$), but now {\it solvable\/}
for any $(I^+,I^-)$-con\-den\-ser~$\mathcal A$ (e.\,g., even
nonclosed or infinite) and any vector~$a$. We have succeeded in this
under the following conditions, which will always be tacitly
assumed.

From now on, in addition to~(\ref{nonzero1}), the following {\bf
standing assumptions} will be always required. The kernel~$\kappa$
is assumed to be consistent, and either
$$
I^-=\varnothing,
$$
or the following three conditions are satisfied:
\begin{equation}
g_{\inf}>0,\label{g}
\end{equation}
\begin{equation}\sup_{x\in A^+,\ y\in
A^-}\,\kappa(x,y)<\infty,\label{bou}
\end{equation}
\begin{equation}|a|:=\sum_{i\in I}\,a_i<\infty.\label{abou}
\end{equation}

{\bf Remark 10.} These assumptions on a kernel are not too
restrictive. In particular, they all are satisfied by the Newtonian,
Riesz, or Green kernels in~$\mathbb R^n$, $n\geqslant2$, provided
the Euclidean distance between $A^+$ and $A^-$ is nonzero, as well
as by the restriction of the logarithmic kernel in~$\mathbb R^2$ to
an open unit ball.\bigskip

{\bf\large 6. $\mathcal A$-vague and strong cluster sets of
minimizing nets}\nopagebreak\medskip

To formulate the results obtained, we shall need the following
notation.

{\bf 6.1.} Denote by $\mathbb M(\mathcal A,a,g)$ the class of all
$(\mu_t)_{t\in T}\subset\mathcal E^0(\mathcal A,a,g)$ such that
\begin{equation}
\lim_{t\in T}\,\|\mu_t\|^2=\|\mathcal E(\mathcal A,a,g)\|^2.
\label{min}
\end{equation}
This class is not empty, which is clear from~(\ref{nonzero1}) in
view of Corollary~3.

Let $\mathcal M(\mathcal A,a,g)$ (respectively, $\mathcal
M'(\mathcal A,a,g)$) consist of all limit points of the nets
$(\mu_t)_{t\in T}\in\mathbb M(\mathcal A,a,g)$ in the $\mathcal
A$-vague topology of the space~$\mathfrak M(\,\overline{\mathcal
A}\,)$ (respectively, in the strong topology of the semimetric space
$\mathcal E(\,\overline{\mathcal A}\,)$). Also write
$$
\mathcal E(\mathcal A,\leqslant\!a,g):=\Bigl\{\mu\in\mathcal
E(\mathcal A):\quad\int g\,d\mu^i\leqslant a_i\mbox{ \ for all \ }
i\in I\Bigr\}.
$$

With the preceding notation and under our standing assumptions
(see~Sec.~5.2), there holds the following lemma, to be proved in
Sec.~12 below.\medskip

{\bf Lemma 8.} {\it Given $(\mu_t)_{t\in T}\in\mathbb M(\mathcal
A,a,g)$, there exist its $\mathcal A$-vague cluster points; hence,
$\mathcal M(\mathcal A,a,g)$ is nonempty. Moreover,
\begin{equation}
\mathcal M(\mathcal A,a,g)\subset\mathcal M'(\mathcal
A,a,g)\cap\mathcal E(\,\overline{\mathcal
A},\leqslant\!a,g).\label{MMprime}
\end{equation}
Furthermore, for every $\chi\in\mathcal M'(\mathcal A,a,g)$,
\begin{equation}
\lim_{t\in T}\,\|\mu_t-\chi\|^2=0,\label{strongly}
\end{equation}
and hence $\mathcal M'(\mathcal A,a,g)$ forms an equivalence class
in~$\mathcal E(\,\overline{\mathcal A}\,)$.}\medskip

It follows from (\ref{min})\,--\,(\ref{strongly}) that
$$
\|\zeta\|^2=\|\mathcal E(\mathcal A,a,g)\|^2\quad\mbox{for all \
}\zeta\in\mathcal M(\mathcal A,a,g).
$$
Also observe that, if $\mathcal A=\mathcal K$ is compact, then
moreover $\mathcal M(\mathcal K,a,g)\subset\mathfrak M(\mathcal
K,a,g)$, which together with the preceding relation proves the
following assertion.\medskip

{\bf Corollary 6.} {\it If $\mathcal A=\mathcal K$ is compact, then
the $\mathcal E(\mathcal K,a,g)$-problem is solvable. Actually,}
\begin{equation}
\mathcal S(\mathcal K,a,g)=\mathcal M(\mathcal K,a,g).\label{S}
\end{equation}

{\bf 6.2.} When approaching $\mathcal A$ by the increasing
family~$\{\mathcal K\}_\mathcal A$ of the compact condensers
$\mathcal K\prec\mathcal A$, we shall always suppose all
those~$\mathcal K$ to be of capacity nonzero. This involves no loss
of generality, which is clear from~(\ref{nonzero1}) and Lemma~5.

Then Corollary 6 enables us to introduce the (nonempty) class
$\mathbb M_0(\mathcal A,a,g)$ of all nets $(\lambda_\mathcal
K)_{\mathcal K\in\{\mathcal K\}_\mathcal A}$, where
$\lambda_\mathcal K\in\mathcal S(\mathcal K,a,g)$ is arbitrarily
chosen. Let $\mathcal M_0(\mathcal A,a,g)$ consist of all $\mathcal
A$-vague cluster points of those nets. Since, by Lemma~5,
$$\mathbb M_0(\mathcal A,a,g)\subset\mathbb M(\mathcal A,a,g),$$
application of Lemma~8 yields the following assertion.\medskip

{\bf Corollary 7.} {\it The class $\mathcal M_0(\mathcal A,a,g)$ is
nonempty, and}
$$
\mathcal M_0(\mathcal A,a,g)\subset\mathcal M(\mathcal
A,a,g)\subset\mathcal M'(\mathcal A,a,g).
$$

{\bf Remark 11.} Each of the cluster sets, $\mathcal M_0(\mathcal
A,a,g)$, $\mathcal M(\mathcal A,a,g)$ and $\mathcal M'(\mathcal
A,a,g)$, plays an important role in our study. However, if $\kappa$
is additionally assumed to be strictly positive definite (hence,
perfect), while $\overline{A_i}$, $i\in I$, are mutually disjoint,
then all these three classes coincide and consist of just one
element.\medskip

{\bf 6.3.} Also the following notation will be required. Given
$\chi\in\mathcal M'(\mathcal A,a,g)$, write
$$
\mathcal M'_\mathcal E(\mathcal A,a,g):=\bigl[R\chi\bigr]_\mathcal
E\,.
$$
This equivalence class does not depend on the choice of $\chi$,
which is clear from Lemma~8. Lemma~8 also yields that, for any
$(\mu_t)_{t\in T}\in\mathbb M(\mathcal A,a,g)$ and any
$\nu\in\mathcal M'_\mathcal E(\mathcal A,a,g)$, $R\mu_t\to\nu$ in
the strong topology of the pre-Hilbert space~$\mathcal E$.\bigskip

{\bf\large 7. Extremal problems dual to the main
minimum-pro\-blem}\nopagebreak\medskip

Throughout Sec.~7, as usual, we are keeping all our standing
assumptions, stated in~Sec.~5.2.\medskip

{\bf 7.1.} A proposition $R(x)$ involving a variable point
$x\in\mathbf X$ is said to subsist {\it nearly
everywhere\/}~(n.\,e.) in~$E$, where $E$ is a given subset
of~$\mathbf X$, if the set of all $x\in E$ for which $R(x)$ fails to
hold is of interior capacity zero. See,~e.\,g.,~\cite{F1}.

If $C(E)>0$ and $f$ is a universally measurable function bounded
from below nearly everywhere in~$E$, write
$$
"\!\inf_{x\in E}\!"\,f(x):=\sup\,\bigl\{q: \ f(x)\geqslant
q\quad\mbox{n.\,e.~in \ } E\bigr\}.
$$
Then
$$
f(x)\geqslant"\!\inf_{x\in E}\!"\,f(x)\quad\mbox{n.\,e.~in \ } E.
$$
This follows immediately from the fact, to be often used in what
follows, that the union of a sequence of sets $U_n\cap E$ with
$C(U_n\cap E)=0$ is of interior capacity zero as well, provided
$U_n$, $n\in\mathbb N$, are universally measurable whereas $E$ is
arbitrary (see the corollary to~Lemma~2.3.5 in~\cite{F1} and the
remark attached to it).\medskip

{\bf 7.2.} Let $\hat{\Gamma}=\hat{\Gamma}(\mathcal A,a,g)$ denote
the class of all Radon measures $\nu\in\mathcal E$ such that there
exist real numbers $c_i(\nu)$, $i\in I$, satisfying the relations
\begin{equation}
\alpha_i a_i\kappa(x,\nu)\geqslant c_i(\nu)g(x)\quad\mbox{n.\,e. in
\ } A_i,\quad i\in I, \label{adm1}
\end{equation}
\begin{equation}
\sum_{i\in I}\,c_i(\nu)\geqslant1. \label{adm2}
\end{equation}

{\bf Remark~12.} Given $\nu\in\hat{\Gamma}$, then the series
in~(\ref{adm2}) must be {\it absolutely convergent}. Indeed, due
to~(\ref{nonzero1}) and Corollary~3, there exists $\mu\in\mathcal
E^0(\mathcal A,a,g)$; then, by~\cite[Lemma~2.3.1]{F1}, the
inequality in~(\ref{adm1}) holds $\mu^i$-almost everywhere. In view
of $\int g\,d\mu^i=a_i$, this gives
$$\kappa(\alpha_i\mu^i,\nu)\geqslant c_i(\nu),\quad i\in I.$$
Since, by Fubini's theorem and Lemma~3, $\sum_{i\in
I}\,\kappa(\alpha_i\mu^i,\nu)$ absolutely converges, the required
conclusion follows.\medskip

We also observe that the class $\hat{\Gamma}(\mathcal A,a,g)$ is
{\it convex\/}, which can easily be seen from the property of sets
of interior capacity zero mentioned just above.

The following assertion, to be proved in Sec.~15 below, holds
true.\medskip

{\bf Theorem 2.} {\it Under the standing assumptions,}
\begin{equation}
\|\hat{\Gamma}(\mathcal A,a,g)\|^2={\rm cap}\,(\mathcal A,a,g).
\label{id}
\end{equation}

If $\|\hat{\Gamma}(\mathcal A,a,g)\|^2<\infty$, we are interested in
the $\hat{\Gamma}(\mathcal A,a,g)$-{\it problem\/} (cf.~Sec.~4.1),
i.\,e., the problem on the existence of
$\hat{\omega}\in\hat{\Gamma}(\mathcal A,a,g)$ with minimal energy
$$
\|\hat{\omega}\|^2=\|\hat{\Gamma}(\mathcal A,a,g)\|^2;
$$
the collection of all those $\hat\omega$ will be denoted by
$\hat{\mathcal G}=\hat{\mathcal G}(\mathcal A,a,g)$.

A minimizing measure $\hat{\omega}$ can be shown to be {\it
unique\/} up to a summand of seminorm zero (and, hence, it is unique
whenever the kernel under consideration is strictly positive
definite). Actually, the following stronger result holds
true.\medskip

{\bf Lemma 9.} {\it If $\hat\omega$ exists, $\hat{\mathcal
G}(\mathcal A,a,g)$ forms an equivalence class in~$\mathcal
E$.}\medskip

{\bf Proof.} Since $\hat{\Gamma}$ is convex, Lemma~4 yields that
$\hat{\mathcal G}$ is contained in an equivalence class in~$\mathcal
E$. To prove that $\hat{\mathcal G}$ actually coincides with that
equivalence class, it suffices to show that, if $\nu$ belongs
to~$\hat{\Gamma}$, then so do all measures equivalent to~$\nu$
in~$\mathcal E$. But this follows at once from the property of sets
of interior capacity zero mentioned in Sec.~7.1 and the fact that
the potentials of any two equivalent in~$\mathcal E$ measures
coincide nearly everywhere in~$\mathbf
X$~\cite[Lemma~3.2.1]{F1}.\hfill\phantom{B}\hfill$\Box$\medskip

{\bf 7.3.} Assume for a moment that ${\rm cap}\,(\mathcal A,a,g)$ is
finite. When combined with~(\ref{def}) and~(\ref{iden}), Theorem~2
shows that the $\hat{\Gamma}(\mathcal A,a,g)$-pro\-blem and, on the
other hand, the $\mathcal E(\mathcal A,a\,{\rm cap}\,\mathcal
A,g)$-pro\-blem have the same infimum, equal to the capacity ${\rm
cap}\,\mathcal A$, and so these two variational problems are {\it
dual}.

But what is surprising is that their infimum, ${\rm cap}\,\mathcal
A$, turns out to be always an actual minimum in the former extremal
problem, while this is not the case for the latter one
(see~Sec.~5.1). In fact, the following statement on the solvability
of the $\hat{\Gamma}(\mathcal A,a,g)$-problem, to be proved in
Sec.~15 below, holds true.\medskip

{\bf Theorem 3.} {\it Under the standing assumptions, if moreover
${\rm cap}\,\mathcal A<\infty$, then the class $\hat{\mathcal
G}(\mathcal A,a,g)$ is nonempty and can be given by the formula
\begin{equation}
\hat{\mathcal G}(\mathcal A,a,g)=\mathcal M'_\mathcal E(\mathcal
A,a\,{\rm cap}\,\mathcal A,g).\label{desc}
\end{equation}
The numbers $c_i(\hat{\omega})$, $i\in I$, satisfying
both\/}~(\ref{adm1}) {\it and\/}~(\ref{adm2}) {\it for
$\hat{\omega}\in\hat{\mathcal G}(\mathcal A,a,g)$, are determined
uniquely, do not depend on the choice of~$\hat{\omega}$, and can be
written in either of the forms
\begin{align}
c_i(\hat{\omega})&=\alpha_i\,{\rm cap}\,\mathcal
A^{-1}\kappa(\zeta^i,\zeta),\label{snt1}\\[4pt]
c_i(\hat{\omega})&=\alpha_i\,{\rm cap}\,\mathcal A^{-1}\lim_{s\in
S}\,\kappa(\mu_s^i,\mu_s),\label{cnt}
\end{align}
where $\zeta\in\mathcal M(\mathcal A,a\,{\rm cap}\,\mathcal A,g)$
and $(\mu_s)_{s\in S}\in\mathbb M(\mathcal A,a\,{\rm cap}\,\mathcal
A,g)$ are arbitrarily given.}\medskip

The following two assertions, providing additional information
about~$c_i(\hat{\omega})$, $i\in I$, can be obtained directly from
the preceding theorem.\medskip

{\bf Corollary 8.} {\it Given $\hat{\omega}\in\hat{\mathcal
G}(\mathcal A,a,g)$, it follows that}
\begin{equation}
c_i(\hat{\omega})="\!\inf_{x\in A_i}\!"\,\,\frac{\alpha_i
a_i\kappa(x,\hat{\omega})}{g(x)},\quad i\in
I\smallskip.\label{ess''}
\end{equation}

{\bf Corollary 9.} {\it The inequality\/}~(\ref{adm2}) {\it for
$\hat{\omega}\in\hat{\mathcal G}(\mathcal A,a,g)$ is actually an
equality; i.\,e.,} \begin{equation}\sum_{i\in
I}\,c_i(\hat\omega)=1.\label{1''}\end{equation}

{\bf Remark 13.} Assume for a moment that $C(A_j)=0$ for some $j\in
I$. Then ${\rm cap}\,\mathcal A=0$ according to Corollary~4. On the
other hand, the measure $\nu_0=0$ belongs to $\hat{\Gamma}(\mathcal
A,a,g)$ since it satisfies both~(\ref{adm1}) and~(\ref{adm2}) with
$c_i(\nu_0)$, $i\in I$, where
$$c_j(\nu_0)\geqslant1\quad\mbox{and}\quad c_i(\nu_0)=0,\quad i\ne j.$$
This implies that the identity~(\ref{id}) actually holds true in the
degenerate case $C(A_j)=0$ as well, and then $\hat{\mathcal
G}(\mathcal A,a,g)$ consists of all $\nu\in\mathcal E$ of seminorm
zero. What then, however, fails to hold is the statement on the
uniqueness of~$c_i(\hat{\omega})$.\medskip

{\bf 7.4.} Let $\hat{\Gamma}_*(\mathcal A,a,g)$ consist of all
$\nu\in\hat{\Gamma}(\mathcal A,a,g)$ for which the
inequality~(\ref{adm2}) is actually an equality. By arguments
similar to those that have been applied above, one can see that
$\hat{\Gamma}_*(\mathcal A,a,g)$ is convex, and hence all the
solutions to the minimal energy problem over this class form an
equivalence class in~$\mathcal E$. Combining this with Theorems~2,~3
and Corollary~9 leads to the following assertion.\medskip

{\bf Corollary 10.} {\it Under the standing assumptions,
$$
\|\hat{\Gamma}_*(\mathcal A,a,g)\|^2={\rm cap}\,(\mathcal A,a,g).
$$
If moreover ${\rm cap}\,(\mathcal A,a,g)<\infty$, then the
$\hat{\Gamma}_*(\mathcal A,a,g)$-problem is solvable and the class
$\hat{\mathcal G}_*(\mathcal A,a,g)$ of all its solutions is given
by the formula}
$$\hat{\mathcal G}_*(\mathcal A,a,g)=
\mathcal M'_\mathcal E(\mathcal A,a\,{\rm cap}\,\mathcal A,g).
$$

{\bf Remark 14.} Theorem~2 and Corollary~10 (cf.~also Theorem~4 and
Corollary~12 below) provide new equivalent definitions of the
capacity ${\rm cap}\,(\mathcal A,a,g)$. Note that, in contrast to
the initial definition (cf.~Sec.~4.2), no restrictions on the
supports and total masses of measures from the classes
$\hat{\Gamma}(\mathcal A,a,g)$ or $\hat{\Gamma}_*(\mathcal A,a,g)$
have been imposed; the only restriction involves their potentials.
These definitions of the capacity are actually new even in the
simplest case of a finite, compact condenser; compare with~\cite{O}.
They are not only of obvious academic interest, but seem also to be
important for numerical computations.\medskip

{\bf 7.5.} Our next purpose is to formulate an $\mathcal H$-problem
such that it is still dual to the $\mathcal E(\mathcal A,a\,{\rm
cap}\,\mathcal A,g)$-problem and solvable, but now with $\mathcal H$
consisting of measures associated with a condenser.

Let $\Gamma(\mathcal A,a,g)$ consist of all $\mu\in\mathcal
E(\,\overline{\mathcal A}\,)$ for which both the
relations~(\ref{adm1}) and~(\ref{adm2}) hold (with~$\mu$ in place
of~$\nu$). In other words,
$$
\Gamma(\mathcal
A,a,g):=\left\{\mu\in\mathcal E(\,\overline{\mathcal A}\,):\quad
R\mu\in\hat{\Gamma}(\mathcal A,a,g)\right\}.
$$
Observe that the class $\Gamma(\mathcal A,a,g)$ is {\it convex\/}
and
\begin{equation}
\|\Gamma(\mathcal A,a,g)\|^2\geqslant\|\hat{\Gamma}(\mathcal
A,a,g)\|^2. \label{GammahatGamma}
\end{equation}

We proceed to show that the inequality (\ref{GammahatGamma}) is
actually an equality, and that the minimal energy problem, if
considered over the class $\Gamma(\mathcal A,a,g)$, is still
solvable.\medskip

{\bf Theorem 4.} {\it Under the standing assumptions,
\begin{equation}
\|\Gamma(\mathcal A,a,g)\|^2={\rm cap}\,(\mathcal A,a,g).
\label{id2}
\end{equation}
If moreover ${\rm cap}\,(\mathcal A,a,g)<\infty$, then the
$\Gamma(\mathcal A,a,g)$-problem is solvable and the class $\mathcal
G(\mathcal A,a,g)$ of all its solutions~$\omega$ is given by the
formula}
\begin{equation}
\mathcal G(\mathcal A,a,g)=\mathcal M'(\mathcal A,a\,{\rm
cap}\,\mathcal A,g).\label{desc2}
\end{equation}

{\bf Proof.} We can certainly assume ${\rm cap}\,\mathcal A$ to be
finite, for if not, (\ref{id2})~is obtained directly from~(\ref{id})
and~(\ref{GammahatGamma}). Then, according to Lemma~8 with~$a\,{\rm
cap}\,\mathcal A$ instead of~$a$, the class $\mathcal M'(\mathcal
A,a\,{\rm cap}\,\mathcal A,g)$ is nonempty; fix~$\chi$, one of its
elements. It is clear from its definition and the
identity~(\ref{desc}) that $\chi\in\mathcal E(\,\overline{\mathcal
A}\,)$ and $R\chi\in\hat{\mathcal G}(\mathcal A,a,g)$. Hence,
$\chi\in\Gamma(\mathcal A,a,g)$, and therefore
$$\|\hat{\Gamma}(\mathcal A,a,g)\|^2=\|\chi\|^2\geqslant\|\Gamma(\mathcal A,a,g)\|^2.$$
In view of~(\ref{id}) and~(\ref{GammahatGamma}), this
proves~(\ref{id2}) and, as well, the inclusion
$$\mathcal M'(\mathcal A,a\,{\rm cap}\,\mathcal A,g)\subset\mathcal G(\mathcal A,a,g).$$
But the right-hand side of this inclusion is an equivalence class
in~$\mathcal E(\,\overline{\mathcal A}\,)$, which follows from the
convexity of $\Gamma(\mathcal A,a,g)$ and Lemma~4 in the same manner
as in the proof of Lemma~9. Since, by Lemma~8, also the left-hand
side is an equivalence class in~$\mathcal E(\,\overline{\mathcal
A}\,)$, the two sets must actually be
equal.\hfill\phantom{B}\hfill$\Box$\medskip

{\bf Corollary 11.} {\it If $\mathcal A=\mathcal K$ is compact and
${\rm cap}\,(\mathcal K,a,g)<\infty$, then any solution to the
$\mathcal E(\mathcal K,a\,{\rm cap}\,\mathcal K,g)$-problem gives,
as well, a solution to the $\Gamma(\mathcal
K,a,g)$-problem.}\medskip

{\bf Proof.} This follows from~(\ref{desc2}), when combined
with~(\ref{MMprime}) and~(\ref{S}) for~$a\,{\rm cap}\,\mathcal K$ in
place of~$a$.\hfill\phantom{B}\hfill$\Box$\medskip

{\bf Remark 15.} Assume ${\rm cap}\,\mathcal A<\infty$, and fix
$\omega\in\mathcal G(\mathcal A,a,g)$ and
$\hat{\omega}\in\hat{\mathcal G}(\mathcal A,a,g)$. Since,
by~(\ref{desc}) and~(\ref{desc2}),
$\kappa(x,\omega)=\kappa(x,\hat\omega)$ nearly everywhere
in~$\mathbf X$, the numbers~$c_i(\omega)$, $i\in I$,
satisfying~(\ref{adm1}) and~(\ref{adm2}) for~$\nu=\omega$, are
actually equal to~$c_i(\hat\omega)$. This implies that
relations~\mbox{(\ref{snt1})\,--\,(\ref{1''})} do hold, as well,
for~$\omega$ in place of~$\hat\omega$.\medskip

{\bf Remark 16.} Observe that, in all the preceding assertions, we
have not imposed any restrictions on the topology of~$A_i$, $i\in
I$. So, all the $\hat{\Gamma}(\mathcal A,a,g)$-,
$\hat{\Gamma}_*(\mathcal A,a,g)$-, and $\Gamma(\mathcal
A,a,g)$-problems are solvable even for a nonclosed, infinite
condenser~$\mathcal A$.\medskip

{\bf Remark 17.} If $I=\{1\}$ and $g=1$, Theorems~2\,--\,4 and
Corollary~10 can be derived from~\cite{F1}. Moreover, then one can
choose $\gamma\in\mathcal G(\mathcal A,a,g)$ so that
$$\gamma(\mathbf X)=a_1\,C(A_1),$$ and exactly this kind of measures
was called by B.~Fuglede {\it interior capacitary distributions
associated with the set\/}~$A_1$. However, this fact in general can
not be extended to a condenser~$\mathcal A$ consisting more than one
plate; that is, then
$$
\mathcal G(\mathcal A,a,g)\cap\mathcal E(\,\overline{\mathcal
A},a\,{\rm cap}\,\mathcal A,g)=\varnothing,
$$
which is seen from the unsolvability of the $\mathcal E(\mathcal
A,a\,{\rm cap}\,\mathcal A,g)$-problem.\bigskip

{\bf\Large 8. Interior capacitary constants associated with a
condenser}\nopagebreak\medskip

{\bf 8.1.} Throughout Sec.~8, it is always required that ${\rm
cap}\,(\mathcal A,a,g)<\infty$. Due to the uniqueness statement in
Theorem~3, the following notion naturally arises.\medskip

{\bf Definition 6.} The numbers
$$C_i:=C_i(\mathcal A,a,g):=c_i(\hat\omega),\quad i\in I,$$
satisfying both the relations~(\ref{adm1}) and~(\ref{adm2}) for
$\hat{\omega}\in\hat{\mathcal G}(\mathcal A,a,g)$, are said to be
the ({\it interior}) {\it capacitary constants\/} associated with an
$(I^+,I^-)$-condenser~$\mathcal A$.\medskip

{\bf Corollary 12.} {\it The interior capacity ${\rm cap}\,(\mathcal
A,a,g)$ equals the infimum of~$\kappa(\nu,\nu)$, where $\nu$ ranges
over the class of all $\nu\in\mathcal E$} ({\it similarly,}
$\nu\in\mathcal E(\,\overline{\mathcal A}\,)$) {\it such that
$$
\alpha_ia_i\kappa(x,\nu)\geqslant C_i(\mathcal
A,a,g)\,g(x)\quad\mbox{n.\,e.~in \ } A_i,\quad i\in I.
$$
The infimum is attained at any} $\hat\omega\in\hat{\mathcal
G}(\mathcal A,a,g)$ ({\it respectively,} $\omega\in\mathcal
G(\mathcal A,a,g)$), {\it and hence it is an actual
minimum.}\medskip

{\bf Proof.} This follows immediately from Theorems~2\,--\,4 and
Remark~15.\hfill\phantom{B}\hfill$\Box$\medskip

{\bf 8.2.} Some properties of the interior capacitary constants
$C_i(\mathcal A,a,g)$, $i\in I$, have already been provided by
Theorem~3 and Corollaries~8,~9. Also observe that, if $I=\{1\}$,
then certainly $C_1(\mathcal A,a,g)=1$
(cf.~\cite[Th.~4.1]{F1}).\medskip

{\bf Corollary 13.} {\it $C_i(\,\cdot\,,a,g)$, $i\in I$, are
continuous under exhaustion of~$\mathcal A$ by the increasing family
of all compact condensers $\mathcal K\prec\mathcal A$. Namely,}
$$
C_i(\mathcal A,a,g)=\lim_{\mathcal K\uparrow\mathcal
A}\,C_i(\mathcal K,a,g).$$

{\bf Proof.} Under our assumptions, $0<{\rm cap}\,\mathcal K<\infty$
for every $\mathcal K\in\{\mathcal K\}_\mathcal A$, and hence there
exists $\lambda_\mathcal K\in\mathcal S(\mathcal K,a\,{\rm
cap}\,\mathcal K,g)$. Substituting $\lambda_\mathcal K$
into~(\ref{snt1}) yields
\begin{equation}
C_i(\mathcal K,a,g)=\alpha_i\,{\rm cap}\,\mathcal
K^{-1}\,\kappa(\lambda^i_\mathcal K,\lambda_\mathcal K),\quad i\in
I.\label{cntt}
\end{equation}
On the other hand, it follows from Lemma~5 that the net
$$
{\rm cap}\,\mathcal A\,{\rm cap}\,\mathcal K^{-1}\,\lambda_\mathcal
K,\quad\mbox{where \ } \mathcal K\in\{\mathcal K\}_\mathcal A,
$$
belongs to the class $\mathbb M(\mathcal A,a\,{\rm cap}\,\mathcal
A,g)$. Substituting it into~(\ref{cnt}) and then combining the
relation obtained with~(\ref{cntt}), we get the
corollary.\hfill$\Box$\medskip

{\bf Corollary 14.} {\it Assume $C(A_j)=\infty$ for some $j\in I$.
If moreover $g_{\inf}>0$, then}
$$
C_j(\mathcal A,a,g)\leqslant0.
$$

{\bf Proof.} Assume, on the contrary, that $C_j>0$. Given
$\hat{\omega}\in\hat{\mathcal G}(\mathcal A,a,g)$, then
$$
\alpha_j a_j\kappa(x,\hat{\omega})\geqslant
C_j\,g_{\inf}>0\quad\mbox{n.\,e. in \ } A_j,
$$
and therefore, by \cite[Lemma 3.2.2]{F1},
$$
C(A_j)\leqslant
a_j^2\,\|\hat{\omega}\|^2\,C_j^{-2}\,g_{\inf}^{-2}<\infty,
$$
which is a contradiction.\hfill\phantom{B}\hfill$\Box$\medskip

{\bf Remark 18.} Observe that the necessity part of Lemma~7, which
has been proved above with elementary arguments, can also be
obtained as a consequence of Corollary~14. Indeed, if (\ref{lemma6})
were true, then by~Corollary~14 the sum of~$C_i$, where $i$ ranges
over~$I$, would be not greater than~$0$, which is
impossible.\bigskip

{\bf\Large 9. Interior capacitary distributions associated with a
condenser}\nopagebreak\medskip

As always, we are keeping all our standing assumptions, stated
in~Sec.~5.2. Throughout Sec.~9, it is also required that ${\rm
cap}\,\mathcal A<\infty$.

Our next purpose is to introduce a notion of interior capacitary
distributions~$\gamma_\mathcal A$ associated with a
condenser~$\mathcal A$ such that the distributions obtained possess
properties similar to those of interior capacitary distributions
associated with a set. Fuglede's theory of interior capacities of
sets~\cite{F1} serves here as a model case.

{\bf 9.1.} If $\mathcal A=\mathcal K$ is compact, then, as follows
from Theorem~4, Corollary~11 and Remark~15, any
minimizer~$\lambda_\mathcal K$ in the $\mathcal E(\mathcal K,a\,{\rm
cap}\,\mathcal K,g)$-problem has the desired properties, and so
$\gamma_\mathcal K$ might be defined as $$\gamma_\mathcal
K:=\lambda_\mathcal K,\quad\mbox{where}\quad\lambda_\mathcal
K\in\mathcal S(\mathcal K,a\,{\rm cap}\,\mathcal K,g).$$ However, as
is seen from Remark~17, in the noncompact case the desired notion
can not be obtained as just a direct generalization of the
corresponding one from the theory of interior capacities of sets.
Having in mind that, similar to our model case, the required
distributions should give a solution to the $\Gamma(\mathcal
A,a,g)$-problem and be strongly and $\mathcal A$-vaguely continuous
under exhaustion of~$\mathcal A$ by compact condensers, we arrive at
the following definition.\medskip

{\bf Definition 7.} We shall call $\gamma_\mathcal A\in\mathcal
E(\,\overline{\mathcal A}\,)$ an ({\it interior}) {\it capacitary
distribution\/} associated with~$\mathcal A$ if there exists a
subnet $(\mathcal K_s)_{s\in S}$ of $(\mathcal K)_{\mathcal
K\in\{\mathcal K\}_\mathcal A}$ and
$$\lambda_{\mathcal K_s}\in\mathcal S(\mathcal K_s,a\,{\rm
cap}\,\mathcal K_s,g),\quad s\in S,$$ such that $(\lambda_{\mathcal
K_s})_{s\in S}$ converges to~$\gamma_\mathcal A$ in both the
$\mathcal A$-vague and the strong topologies. Let $\mathcal
D(\mathcal A,a,g)$ denote the collection of all
those~$\gamma_\mathcal A$.\medskip

Application of Lemmas 5 and 8 enables us to rewrite the above
definition in the following, apparently weaker, form:
\begin{equation}
\mathcal D(\mathcal A,a,g)=\mathcal M_0(\mathcal A,a\,{\rm
cap}\,\mathcal A,g).\label{D}
\end{equation}

{\bf Theorem 5.} \ {\it $\mathcal D(\mathcal A,a,g)$ is nonempty,
contained in an equivalence class in~$\mathcal
E(\,\overline{\mathcal A}\,)$, and compact in the induced $\mathcal
A$-vague topology. Furthermore,
\begin{equation}
\mathcal D(\mathcal A,a,g)\subset\mathcal G(\mathcal
A,a,g)\cap\mathcal E(\,\overline{\mathcal A},\leqslant\!a\,{\rm
cap}\,\mathcal A,g).\label{gammain}
\end{equation}
Given $\gamma:=\gamma_\mathcal A\in\mathcal D(\mathcal A,a,g)$, then
\begin{equation}
\|\gamma\|^2={\rm cap}\,\mathcal A, \label{5.1}
\end{equation}
\begin{equation}
\alpha_i a_i\kappa(x,\gamma)\geqslant C_i\,g(x)\quad\mbox{n.\,e. in
\ } A_i,\quad i\in I,\label{5.3}
\end{equation}
where $C_i=C_i(\mathcal A,a,g)$, $i\in I$, are the interior
capacitary constants. Actually,
\begin{equation}
C_i=\frac{\alpha_i\kappa(\gamma^i,\gamma)}{{\rm cap}\,\mathcal
A}="\!\inf_{x\in A_i}\!"\,\,\frac{\alpha_i
a_i\kappa(x,\gamma)}{g(x)}\,,\quad i\in I.\label{5.2}
\end{equation}
If $I^-\ne\varnothing$, assume moreover that the kernel
$\kappa(x,y)$ is continuous on $\overline{A^+}\times\overline{A^-}$,
while\/} $\kappa(\cdot,y)\to0$ ({\it as\/}~$y\to\infty$) {\it
uniformly on compact sets. Then, for every $i\in I$,}
\begin{align}
\alpha_i a_i\kappa(x,\gamma)&\leqslant
C_i\,g(x)&&\!\!\!\!\!\!\!\!\!\!\!\!\!\!\!\!\!\!\!\!\!\!\!\!\!\!\!\!\!\!\text{\it
for all \ } x\in
S(\gamma^i),\label{5.4}\\
\intertext{\it and hence} \alpha_i
a_i\kappa(x,\gamma)&=C_i\,g(x)&&\!\!\!\!\!\!\!\!\!\!\!\!\!\!\!\!\!\!\!\!\!\!\!\!\!\!\!\!\!\!\text{\it
n.\,e. in \ } A_i\cap S(\gamma^i).\notag
\end{align}

Also note that $\mathcal D(\mathcal A,a,g)$ is contained in an
equivalence class in~$\mathfrak M(\,\overline{\mathcal A}\,)$
provided the kernel~$\kappa$ is strictly positive definite, and it
consists of a {\it unique\/} element~$\gamma_\mathcal A$ if,
moreover, all $\overline{A_i}$, $i\in I$, are mutually
disjoint.\medskip

{\bf Remark 19.} As is seen from the preceding theorem, the
properties of interior capacitary distributions associated with a
condenser are quite similar to those of interior capacitary
distributions associated with a set (cf.~\cite[Th.~4.1]{F1}). The
only important difference is that the sign~$\leqslant$ in the
inclusion~(\ref{gammain}) in general can not be omitted~--- even for
a finite, closed, noncompact condenser. Cf.~Remark~17.\medskip

{\bf Remark 20.} Like as in the theory of interior capacities of
sets, in general none of the $i$-coordinates of $\gamma_\mathcal A$
is concentrated on~$A_i$ (unless $A_i$ is closed). Indeed, let
$\mathbf X=\mathbb R^n$, $n\geqslant3$, $\kappa(x,y)=|x-y|^{2-n}$,
$g=1$, $I^+=\{1\}$, $I^-=\{2\}$, $a_1=a_2=1$, and let $A_1=\{x:
|x|<r\}$ and $A_2=\{x: |x|>R\}$, where $0<r<R<\infty$. Then it can
be shown that
$$\gamma_\mathcal
A=\gamma_{\overline{\mathcal A}}=\bigl[\theta^+-\theta^-\bigr]\,{\rm
cap}\,\mathcal A,$$ where $\theta^+$ and~$\theta^-$ are obtained by
the uniform distribution of unit mass over the spheres~$S(0,r)$
and~$S(0,R)$, respectively. Hence, $|\gamma_\mathcal
A|(A)=0$.\medskip

{\bf 9.2.} The purpose of this section is to point out
characteristic properties of the interior capacitary distributions
and the interior capacitary constants.\medskip

{\bf Proposition 1.} {\it Assume $\mu\in\mathcal
E(\,\overline{\mathcal A}\,)$ has the properties
$$
\|\mu\|^2={\rm cap}\,(\mathcal A,a,g),
$$
$$
\alpha_i
a_i\kappa(x,\mu)\geqslant\frac{\alpha_i\kappa(\mu^i,\mu)}{{\rm
cap}\,\mathcal A}\,g(x)\quad\mbox{\it n.\,e. in \ } A_i,\quad i\in
I.
$$
Then $\mu$ is equivalent in $\mathcal E(\,\overline{\mathcal A}\,)$
to every $\gamma_\mathcal A\in\mathcal D(\mathcal A,a,g)$ and, for
all $i\in I$,}
$$
C_i(\mathcal A,a,g)=\frac{\alpha_i\kappa(\mu^i,\mu)}{{\rm
cap}\,\mathcal A}="\!\inf_{x\in A_i}\!"\,\,\frac{\alpha_i
a_i\kappa(x,\mu)}{g(x)}\,.
$$

Actually, there holds the following stronger result, to be proved in
Sec.~17 below.\medskip

{\bf Proposition 2.} {\it Let $\nu\in\mathcal E(\,\overline{\mathcal
A}\,)$ and $\tau_i\in\mathbb R$, $i\in I$, satisfy the relations
\begin{equation}
\alpha_i a_i\kappa(x,\nu)\geqslant\tau_i\,g(x)\quad\mbox{\it n.\,e.
in \ } A_i,\quad i\in I,\label{desc3'}
\end{equation}
\begin{equation}
\sum_{i\in I}\,\tau_i=\frac{{\rm cap}\,\mathcal A+\|\nu\|^2}{2\,{\rm
cap}\,\mathcal A}\,.\label{adm2''}
\end{equation}
Then $\nu$ is equivalent in $\mathcal E(\,\overline{\mathcal A}\,)$
to every $\gamma_\mathcal A\in\mathcal D(\mathcal A,a,g)$ and, for
all $i\in I$,}
\begin{equation}
\tau_i=C_i(\mathcal A,a,g)="\!\inf_{x\in A_i}\!"\,\,\frac{\alpha_i
a_i\kappa(x,\nu)}{g(x)}\,.\label{un}
\end{equation}

Thus, under the conditions of Proposition~1 or~2, if moreover
$\kappa$ is strictly positive definite and all $\overline{A_i}$,
$i\in I$, are mutually disjoint, then the measure under
consideration is actually the (unique) interior capacitary
distribution~$\gamma_\mathcal A$.\bigskip

{\bf\Large 10. On continuity of the interior capacities, capacitary
distributions, and capacitary constants}\nopagebreak\medskip

{\bf 10.1.} Given $\mathcal A_n=(A_i^n)_{i\in I}$, $n\in\mathbb N$,
and $\mathcal A$ in $\mathfrak C=\mathfrak C(I^+,I^-)$, we write
$\mathcal A_n\uparrow\mathcal A$ if $\mathcal A_n\prec\mathcal
A_{n+1}$ for all~$n$ and
$$A_i=\bigcup_{n\in\mathbb N}\,A_i^n,\quad i\in I.$$

Following \cite[Chap.~1, \S\,9]{B1}, we call a locally compact space
{\it countable at infinity\/} if it can be written as a countable
union of compact sets.\medskip

{\bf Theorem 6.} {\it Suppose that either $g_{\inf}>0$ or the space
$\mathbf X$ is countable at infinity. If $\mathcal A_n$,
$n\in\mathbb N$, are universally measurable and $\mathcal
A_n\uparrow\mathcal A$, then
\begin{equation}
{\rm cap}\,(\mathcal A,a,g)=\lim_{n\in\mathbb N}\,{\rm
cap}\,(\mathcal A_n,a,g). \label{6.1}
\end{equation}
Assume moreover ${\rm cap}\,(\mathcal A,a,g)$ to be finite, and let
$\gamma_n:=\gamma_{\mathcal A_n}$, $n\in\mathbb N$, denote an
interior capacitary distribution associated with~$\mathcal A_n$. If
$\gamma$ is an $\mathcal A$-vague limit point
of\/}~$(\gamma_n)_{n\in\mathbb N}$ ({\it such a $\gamma$ exists\/}),
{\it then $\gamma$ is actually an interior capacitary distribution
associated with the condenser~$\mathcal A$, and
$$ \lim_{n\in\mathbb
N}\,\|\gamma_n-\gamma\|^2=0.
$$
Furthermore,}
\begin{equation}
C_i(\mathcal A,a,g)=\lim_{n\in\mathbb N}\,C_i(\mathcal
A_n,a,g),\quad i\in I\smallskip.\label{6.3}
\end{equation}

Thus, if $\kappa$ is additionally assumed to be strictly positive
definite (hence, perfect) and all $\overline{A_i}$, $i\in I$, are
mutually disjoint, then the (unique) interior capacitary
distribution associated with~$\mathcal A_n$ converges both $\mathcal
A$-vaguely and strongly to the (unique) interior capacitary
distribution associated with~$\mathcal A$.\medskip

{\bf Remark 21.} Theorem 6 remains true if $(\mathcal
A_n)_{n\in\mathbb N}$ is replaced by the increasing ordered family
of all compact condensers~$\mathcal K$ such that $\mathcal
K\prec\mathcal A$. Moreover, then the assumption that either
$g_{\inf}>0$ or $\mathbf X$ is countable at infinity can be omitted.
Cf., e.\,g., Lemma~5 and Corollary~13.\medskip

{\bf Remark 22.} If $I=\{1\}$ and $g=1$, Theorem~6 has been proved
in~\cite[Th.~4.2]{F1}.\medskip

{\bf 10.2.} The remainder of the article is devoted to proving the
results formulated in Sec.~\mbox{6\,--\,10} and is organized as
follows. Theorems~2,~3,~5, and~6 are proved in~Sec.~15,~16, and~18.
Their proofs utilize a description of the potentials of measures
from the classes $\mathcal M'(\mathcal A,a,g)$ and $\mathcal
M_0(\mathcal A,a,g)$, to be given in Sec.~13 and~14 by Lemmas~12
and~13. In turn, Lemmas~12 and~13 use a theorem on the strong
completeness of proper subspaces of~$\mathcal E(\,\overline{\mathcal
A}\,)$, which is a subject of Sec.~11.\bigskip

{\bf\Large 11. On the strong completeness}\nopagebreak\medskip

{\bf 11.1.} Keeping all our standing assumptions on~$\kappa$, $g$,
$a$, and~$\mathcal A$, stated in~Sec.~5.2, we consider $\mathcal
E(\,\overline{\mathcal A},\leqslant\!a,g)$ to be a topological
subspace of the semimetric spa\-ce~$\mathcal E(\,\overline{\mathcal
A}\,)$; the induced topology is likewise called the {\it strong\/}
topology.\medskip

{\bf Theorem 7.} {\it Suppose $\mathcal A$ is closed. Then the
semimetric space $\mathcal E(\mathcal A,\leqslant\!a,g)$ is
complete. In more detail, if $(\mu_s)_{s\in S}\subset\mathcal
E(\mathcal A,\leqslant\!a,g)$ is a strong Cauchy net and $\mu$~is
its $\mathcal A$-vague cluster point\/} ({\it such a $\mu$
exists\/}), {\it then $\mu\in\mathcal E(\mathcal A,\leqslant\!a,g)$
and
\begin{equation}
\lim_{s\in S}\,\|\mu_s-\mu\|^2=0.\label{str}
\end{equation}
Assume, in addition, that the kernel is strictly positive definite
and all $A_i$, $i\in I$, are mutually disjoint. If moreover
$(\mu_s)_{s\in S}\subset\mathcal E(\mathcal A,\leqslant\!a,g)$
converges strongly to $\mu_0\in\mathcal E(\mathcal A)$, then
actually $\mu_0\in\mathcal E(\mathcal A,\leqslant\!a,g)$ and
$\mu_s\to\mu_0$ $\mathcal A$-vaguely}.\medskip

{\bf Remark 23.} This theorem is certainly of independent interest
since, according to the well-known counterexample by
H.~Cartan~\cite{Car}, the pre-Hilbert space~$\mathcal E$ is strongly
incomplete even for the Newtonian kernel $|x-y|^{2-n}$ in~$\mathbb
R^n$, $n\geqslant3$.\medskip

{\bf Remark 24.} Assume the kernel is strictly positive definite
(hence, perfect). If moreover $I^-=\varnothing$, then Theorem~7
remains valid for $\mathcal E(\mathcal A)$ in place of $\mathcal
E(\mathcal A,\leqslant\!a,g)$ (cf.~Theorem~1).  A question still
unanswered is {\bf whether this is the case if $I^+$ and $I^-$ are
both nonempty\/}. We can however show that this is really so for the
Riesz kernels $|x-y|^{\alpha-n}$, $0<\alpha<n$, in~$\mathbb R^n$,
$n\geqslant2$ (cf.~\cite[Th.~1]{Z2}). The proof utilizes Deny's
theorem~\cite{D1} stating that, for the Riesz kernels, $\mathcal E$
can be completed with making use of distributions of finite
energy.\medskip

{\bf 11.2.} We start by auxiliary assertions to be used in the proof
of Theorem~7.\medskip

{\bf Lemma 10.} \ {\it $\mathcal E(\mathcal A,\leqslant\!a,g)$ is
$\mathcal A$-vaguely bounded.}\medskip

{\bf Proof.} Fix $i\in I$, and let a compact set $K\subset A_i$ be
given. Since $g$ is positive and continuous, the inequalities
$$
a_i\geqslant\int g\,d\mu^i\geqslant\mu^i(K)\,\min_{x\in K}\
g(x),\quad\mbox{where \ }\mu\in\mathcal E(\mathcal
A,\leqslant\!a,g),
$$
yield
$$
\sup_{\mu\in\mathcal E(\mathcal A,\leqslant a,g)}\,\mu^i(K)<\infty,
$$
and the lemma follows.\hfill\phantom{B}\hfill$\Box$\medskip

{\bf Lemma 11.} {\it Suppose $\mathcal A$ is closed. If a net
$(\mu_s)_{s\in S}\subset\mathcal E(\mathcal A,\leqslant\!a,g)$ is
strongly bounded, then its $\mathcal A$-vague cluster set is
nonempty and contained in $\mathcal E(\mathcal
A,\leqslant\!a,g)$.}\medskip

{\bf Proof.} We begin by showing that the nets $(R\mu_s^+)_{s\in S}$
and $(R\mu_s^-)_{s\in S}$ are strongly bounded as well, i.\,e.,
\begin{equation}
\sup_{s\in S}\,\|R\mu_s^\pm\|^2<\infty,\quad i\in I.  \label{7.1}
\end{equation}
Of course, this needs to be proved only when $I^-\ne\varnothing$;
then, in accordance with the standing assumptions, all the
relations~(\ref{g}), (\ref{bou}), and~(\ref{abou}) hold. Since
\begin{equation}
\int g\,d\mu_s^i\leqslant a_i,\quad i\in I, \label{oh}
\end{equation}
(\ref{g}) implies
\begin{equation}
\sup_{s\in S}\,\mu_s^i(\mathbf X)\leqslant
a_i\,g_{\inf}^{-1}<\infty,\quad i\in I. \label{7.3}
\end{equation}
Hence, by (\ref{abou}), $$\sup_{s\in S}\,R\mu_s^\pm(\mathbf
X)\leqslant |a|\,g_{\inf}^{-1}<\infty.
$$
When combined with~(\ref{bou}), this shows that
$\kappa(R\mu^+_s,R\mu^-_s)$ is bounded from above on~$S$, and hence
so do $\|R\mu^+_s\|^2$ and $\|R\mu^-_s\|^2$.

Moreover, according to Lemmas~1 and~10, there exists an $\mathcal
A$-vague cluster point~$\mu$ of the net $(\mu_s)_{s\in S}$. Denoting
by $(\mu_d)_{d\in D}$ a subnet of $(\mu_s)_{s\in S}$ such that
\begin{equation}
\mu_d^i\to\mu^i\quad\mbox{vaguely for all \ } i\in I,\label{vagi}
\end{equation}
we get from Lemma 2
\begin{equation}
R\mu_d^\pm\to R\mu^\pm\quad\mbox{vaguely}.\label{vagpm}
\end{equation}
It remains to show that $R\mu^+$ and $R\mu^-$ are both of finite
energy and that (\ref{oh}) holds true for~$\mu^i$ in place
of~$\mu_s^i$. To this end, recall that, if $\mathbf Y$~is a locally
compact Hausdorff space and $\psi$~is a lower semicontinuous
function on~$\mathbf Y$ such that $\psi\geqslant0$ unless its
support is compact, then the map
$$
\nu\mapsto\int\psi\,d\nu,\quad\nu\in\mathfrak M^+(\mathbf Y),
$$
is lower semicontinuous in the induced vague topology (see,
e.\,g.,~\cite{F1}). Applying this to $\mathbf Y=\mathbf
X\times\mathbf X$, $\psi=\kappa$ and, subsequently, $\mathbf Y=A_i$,
$\psi=g|_{A_i}$ and using~(\ref{7.1}), (\ref{vagpm}) and,
respectively, (\ref{oh}) and (\ref{vagi}), we arrive at the required
assertions.\hfill\phantom{B}\hfill$\Box$\medskip

{\bf Corollary 15.} {\it Assume a net $(\mu_s)_{s\in
S}\subset\mathcal E(\mathcal A,\leqslant\!a,g)$ is strongly bounded.
Then for every $i\in I$, $\|\mu^i_s\|^2$ and $\kappa(\mu^i_s,\mu_s)$
are bounded on~$S$.}\medskip

{\bf Proof.} In view of~(\ref{7.1}), the required relation
\begin{equation}
\sup_{s\in S}\,\|\mu_s^i\|^2<\infty\quad\mbox{for all \ } i\in I,
\label{7.1i}
\end{equation}
will be proved once we show that
\begin{equation}\sum_{i,j\in
I^\pm}\,\kappa(\mu_s^i,\mu_s^j)\geqslant
C>-\infty,\label{bdbl}\end{equation} where $C$ is independent
of~$s$. Since (\ref{bdbl}) is obvious when $\kappa\geqslant0$, one
can assume $\mathbf X$ to be compact. Then $\kappa$, being lower
semicontinuous, is bounded from below on~$\mathbf X$ (say by~$-c$,
where $c>0$), while $\mathcal A$ and, hence, $|a|$ are finite.
Furthermore, then $g_{\inf}>0$; therefore, (\ref{7.3}) holds true.
This implies that
$$\kappa(\mu_s^i,\mu_s^j)\geqslant-a_ia_jg_{\inf}^{-2}\,c\quad\mbox{for all \ }i,\,j\in
I,$$ and (\ref{bdbl}) follows.

The above arguments also show that $\kappa(\mu_s^i,R\mu_s^+)$ and
$\kappa(\mu_s^i,R\mu_s^-)$, where $i\in I$ is given, are both
bounded from below on~$S$. Since these functions of~$s$ are bounded
from above as well, which is clear from~(\ref{7.1}) and~(\ref{7.1i})
by the Cauchy-Schwarz inequality, the required boundedness of
$\kappa(\mu_s^i,\mu_s)$ on~$S$ follows.\hfill$\Box$\medskip

{\bf 11.3. Proof of Theorem 7.} Suppose $\mathcal A$ is closed, and
let $(\mu_s)_{s\in S}$ be a strong Cauchy net in $\mathcal
E(\mathcal A,\leqslant\!a,g)$. Since such a net converges strongly
to every its strong cluster point, $(\mu_s)_{s\in S}$ can certainly
be assumed to be strongly bounded. Then, by Lemma~11, there exists
an $\mathcal A$-vague cluster point~$\mu$ of~$(\mu_s)_{s\in S}$, and
\begin{equation}
\mu\in\mathcal E(\mathcal A,\leqslant\!a,g). \label{leqslant}
\end{equation}
We next proceed to verify (\ref{str}). Of course, there is no loss
of generality in assuming $(\mu_s)_{s\in S}$ to converge $\mathcal
A$-vaguely to~$\mu$. Then, by Lemma~2, $(R\mu^+_s)_{s\in S}$ and
$(R\mu^-_s)_{s\in S}$ converge vaguely to~$R\mu^+$ and~$R\mu^-$,
respectively. Since, by~(\ref{7.1}), these nets are strongly bounded
in~$\mathcal E^+$, the property~$(CW)$ (see~Sec.~2) shows that they
approach~$R\mu^+$ and~$R\mu^-$, respectively, in the weak topology
as well, and so
$$
R\mu_s\to R\mu\quad\mbox{weakly}.
$$
This gives
$$
\|\mu_s-\mu\|^2=\|R\mu_s-R\mu\|^2=\lim_{l\in
S}\,\kappa(R\mu_s-R\mu,R\mu_s-R\mu_l),
$$
and hence, by the Cauchy-Schwarz inequality,
$$\|\mu_s-\mu\|^2\leqslant
\|\mu_s-\mu\|\,\liminf_{l\in S}\,\|\mu_s-\mu_l\|,
$$
which proves (\ref{str}) as required, because $\|\mu_s-\mu_l\|$
becomes arbitrarily small when $s,\,l\in S$ are both large enough.

Suppose now that $\kappa$ is strictly positive definite, while all
$A_i$, $i\in I$, are mutually disjoint, and let the net
$(\mu_s)_{s\in S}$ converge strongly to some $\mu_0\in\mathcal
E(\mathcal A)$. Given an $\mathcal A$-vague limit point~$\mu$
of~$(\mu_s)_{s\in S}$, then we conclude from~(\ref{str}) that
$\|\mu_0-\mu\|=0$, hence $\mu_0\cong\mu$ since $\kappa$ is strictly
positive definite, and finally $\mu_0\equiv\mu$ because $A_i$, $i\in
I$, are mutually disjoint. In view of~(\ref{leqslant}), this means
that $\mu_0\in\mathcal E(\mathcal A,\leqslant\!a,g)$, which is a
part of the desired conclusion. Moreover, $\mu_0$ has thus been
shown to be identical to any $\mathcal A$-vague cluster point
of~$(\mu_s)_{s\in S}$. Since the $\mathcal A$-vague topology is
Hausdorff, this implies that $\mu_0$ is actually the $\mathcal
A$-vague limit of~$(\mu_s)_{s\in S}$ (cf.~\cite[Chap.~I, \S~9,
n$^\circ$\,1, cor.]{B1}), which completes the
proof.\hfill\phantom{B}\hfill$\Box$\bigskip

{\bf\Large 12. Proof of Lemma 8}\nopagebreak\medskip

Fix any $(\mu_s)_{s\in S}$ and $(\nu_t)_{t\in T}$ in $\mathbb
M(\mathcal A,a,g)$. It follows by standard arguments that
\begin{equation}
\lim_{(s,t)\in S\times T}\,\|\mu_s-\nu_t\|^2=0, \label{fund}
\end{equation}
where $S\times T$ is the directed product of the directed sets~$S$
and~$T$ (see, e.\,g.,~\cite[Chap.~2,~\S~3]{K}). Indeed, by the
convexity of the class $\mathcal E(\mathcal A,a,g)$,
$$
2\,\|\mathcal E(\mathcal A,a,g)\|
\leqslant{\|\mu_s+\nu_t\|}\leqslant\|\mu_s\|+\|\nu_t\|,
$$
and hence, by (\ref{min}),
$$
\lim_{(s,t)\in S\times T}\,\|\mu_s+\nu_t\|^2=4\,\|\mathcal
E(\mathcal A,a,g)\|^2.
$$
Then the parallelogram identity gives~(\ref{fund}) as claimed.

Relation~(\ref{fund}) implies that $(\mu_s)_{s\in S}$ is strongly
fundamental. Therefore Theorem~7 shows that there exists an
$\mathcal A$-vague cluster point~$\mu_0$ of~$(\mu_s)_{s\in S}$, and
moreover $\mu_0\in\mathcal E(\,\overline{\mathcal
A},\leqslant\!a,g)$ and $\mu_s\to\mu_0$ strongly. This means that
$\mathcal M(\mathcal A,a,g)$ and $\mathcal M'(\mathcal A,a,g)$ are
both nonempty and satisfy the inclusion~(\ref{MMprime}).

What is left is to prove that $\mu_s\to\chi$ strongly, where
$\chi\in\mathcal M'(\mathcal A,a,g)$ is given. But then one can
choose a net in~$\mathbb M(\mathcal A,a,g)$, say $(\nu_t)_{t\in T}$,
convergent to~$\chi$ strongly, and repeated application
of~(\ref{fund}) gives immediately the desired
conclusion.\hfill\phantom{B}\hfill$\Box$\bigskip

{\bf\Large 13. Potentials of strong cluster points of minimizing
nets}\nopagebreak\medskip

{\bf 13.1.} The aim of this section is to provide a description of
the potentials of measures from the class $\mathcal M'(\mathcal
A,a,g)$. As usual, we are keeping all our standing assumptions,
stated in Sec.~5.2.\medskip

{\bf Lemma 12.} {\it There exist $\eta_i\in\mathbb R$, $i\in I$,
such that, for every $\chi\in\mathcal M'(\mathcal A,a,g)$,
\begin{equation}
\alpha_ia_i\kappa(x,\chi) \geqslant\alpha_i\eta_i
g(x)\quad\mbox{n.\,e. in \ } A_i,\quad i\in I, \label{1.1}
\end{equation}
\begin{equation}
\sum_{i\in I}\,\alpha_i\eta_i=\|\mathcal E(\mathcal A,a,g)\|^2.
\label{1.2}
\end{equation}
These $\eta_i$, $i\in I$, are determined uniquely and given by
either of the formulas
\begin{align}
\eta_i&=\kappa(\zeta^i,\zeta),\label{1.3'}\\[2pt]
\eta_i&=\lim_{s\in S}\,\kappa(\mu_s^i,\mu_s),\label{1.3''}
\end{align}
where $\zeta\in\mathcal M(\mathcal A,a,g)$ and $(\mu_s)_{s\in
S}\in\mathbb M(\mathcal A,a,g)$ are arbitrarily chosen.}\medskip

{\bf Proof.} Throughout the proof, we shall assume every net of the
class $\mathbb M(\mathcal A,a,g)$ to be strongly bounded, which
certainly involves no loss of generality.

Fix $\zeta\in\mathcal M(\mathcal A,a,g)$ and choose $(\mu_t)_{t\in
T}\in\mathbb M(\mathcal A,a,g)$ that converges $\mathcal A$-vaguely
to~$\zeta$. We begin by showing that
\begin{equation}
\kappa(\zeta^i,\zeta)=\lim_{t\in T}\,\kappa(\mu_t^i,\mu_t),\quad
i\in I.\label{proof}
\end{equation}
Since, by Corollary~15, $\|\mu_t^i\|$ is bounded from above on~$T$
(say by~$M_1$), while $\mu_t^i\to\zeta^i$ vaguely, the
property~$(CW)$ yields that $\mu_t^i$ approaches~$\zeta^i$ also
weakly. Hence, for every $\varepsilon>0$,
$$
|\kappa(\zeta^i-\mu_t^i,\zeta)|<\varepsilon
$$
whenever $t\in T$ is large enough. Furthermore, by the
Cauchy-Schwarz inequality,
$$
|\kappa(\mu_t^i,\zeta)-\kappa(\mu_t^i,\mu_t)|=|\kappa(\mu_t^i,R\zeta-R\mu_t)|\leqslant
M_1\|\zeta-\mu_t\|,\quad t\in T.
$$
Since, by Lemma 8, $\mu_t\to\zeta$ strongly, the last two relations
combined give~(\ref{proof}).

We next proceed to show that $\eta_i$, $i\in I$, defined by
(\ref{1.3'}), satisfy both~(\ref{1.1}) and~(\ref{1.2}), where
$\chi\in\mathcal M'(\mathcal A,a,g)$ is given. Since (\ref{1.2}) is
obtained directly from
$$\sum_{i\in I}\,\alpha_i\kappa(\zeta^i,\zeta)=\|\zeta\|^2=\|\mathcal E(\mathcal A,a,g)\|^2,$$
suppose,
contrary to~(\ref{1.1}), that there exist $j\in I$ and a set
$E_j\subset A_j$ of interior capacity nonzero such that
\begin{equation}
\alpha_ja_j\kappa(x,\chi)<\alpha_j\eta_jg(x)\quad\mbox{for all \ }
x\in E_j. \label{3.2} \end{equation} Then one can choose
$\nu\in\mathcal E^+$ with compact support so that $S(\nu)\subset
E_j$ and
$$\int g\,d\nu=a_j.$$
Integrating the inequality in (\ref{3.2}) with respect to~$\nu$
gives
\begin{equation}
\alpha_j\,\bigl[\kappa(\chi,\nu)-\eta_j\bigr]<0. \label{3.3}
\end{equation}

To get a contradiction, for every $\tau\in(0,1]$ write
$$
\tilde{\mu}^i_t:=\left\{
\begin{array}{cl} \mu^j_t-\tau\bigl(\mu_t^j-\nu\bigr) & \mbox{if \ }
i=j,\\[4pt]
\mu^i_t & \mbox{otherwise}.\\ \end{array} \right.
$$
Clearly,
$$
\tilde{\mu}_t:=\sum_{i\in I}\,\alpha_i\tilde{\mu}^i_t\in\mathcal
E^0(\mathcal A,a,g),\quad t\in T,
$$
and consequently
\begin{equation} \|\mathcal
E(\mathcal A,a,g)\|^2\leqslant \|\tilde{\mu}_t\|^2=\|\mu_t\|^2
-2\alpha_j\tau\,\kappa(\mu_t,\mu_t^j-\nu)+\tau^2\|\mu_t^j-\nu\|^2.
\label{3.4}
\end{equation}
The coefficient of~$\tau^2$ is bounded from above on~$T$ (say
by~$M_0$), while by Lemma~8
$$
\lim_{t\in T}\,\|\mu_t-\chi\|^2=0.
$$
Combining (\ref{1.3'}), (\ref{proof}) and substituting the result
obtained into~(\ref{3.4}) therefore gives
$$
0\leqslant
M_0\tau^2+2\alpha_j\tau\,\bigl[\kappa(\chi,\nu)-\eta_j\bigr].
$$
By letting here $\tau$ tend to~$0$, we arrive at a contradiction to
(\ref{3.3}), which proves~(\ref{1.1}).

To prove the statement on uniqueness, consider some other $\eta'_i$,
$i\in I$, satisfying both~(\ref{1.1}) and~(\ref{1.2}). Then they are
necessarily finite, and for every~$i$,
\begin{equation}
\alpha_ia_i\kappa(x,\chi)\geqslant\max\bigl\{\alpha_i\eta_i,\,
\alpha_i\eta'_i\bigr\}\,g(x)\quad\mbox{n.\,e. in \ } A_i,\label{4.1}
\end{equation}
which follows from the property of sets of interior capacity zero
mentioned in~Sec.~7.1. Since $\mu_t^i$ is concentrated on~$A_i$ and
has finite energy and compact support, application of~\cite[
Lemma~2.3.1]{F1} shows that the inequality in~(\ref{4.1}) holds
$\mu_t^i$-almost everywhere in~$\mathbf X$. Integrating it with
respect to~$\mu_t^i$ and then summing up over all $i\in I$, in view
of $\int g\,d\mu_t^i=a_i$ we have
$$
\kappa(\mu_t,\chi)\geqslant\sum_{i\in I}\,
\max\bigl\{\alpha_i\eta_i,\,\alpha_i\eta'_i\bigr\},\quad t\in T.
$$
Passing here to the limit as $t$ ranges over $T$, we get
$$
\|\chi\|^2=\lim_{t\in T}\kappa(\mu_t,\chi)\geqslant \sum_{i\in
I}\max\bigl\{\alpha_i\eta_i,\,\alpha_i\eta'_i\bigr\}\geqslant\sum_{i\in
I}\alpha_i\eta_i=\|\mathcal E(\mathcal A,a,g)\|^2,
$$
and hence
$$
\max\bigl\{\alpha_i\eta_i,\,\alpha_i\eta'_i\bigr\}
=\alpha_i\eta_i,\quad i\in I,
$$
for the extreme left and right parts of the above chain of
inequalities are equal. Applying the same arguments again, but with
the roles of $\eta_i$ and $\eta'_i$ reversed, we conclude that
$\eta_i=\eta'_i$ for all $i\in I$, as claimed.

It remains to show that $\eta_i$ can be written in the
form~(\ref{1.3''}), where $(\mu_s)_{s\in S}\in\mathbb M(\mathcal
A,a,g)$ is given. By Corollary~15, for every $i\in I$,
$\kappa(\mu_s^i,\mu_s)$ is bounded on~$S$. Fix $j\in I$ and choose a
cluster point~$\eta_j^0$ of $\bigl\{\kappa(\mu_s^j,\mu_s):\ s\in
S\bigr\}$; then, in view of Lemmas~1 and~10, one can select an
$\mathcal A$-vaguely convergent subnet $(\mu_d)_{d\in D}$
of~$(\mu_s)_{s\in S}$ such that
$$\eta_j^0=\lim_{d\in D}\,\kappa(\mu_d^j,\mu_d).$$
But what has been proved just above implies immediately that
$\eta_j^0=\eta_j$. Since this means that any cluster point of the
net $\kappa(\mu_s^j,\mu_s)$, $s\in S$, coincides with~$\eta_j$,
(\ref{1.3''}) follows.\hfill\phantom{B}\hfill$\Box$\medskip

{\bf 13.2.} In what follows, $\eta_i=:\eta_i(\mathcal A,a,g)$, $i\in
I$, will always denote the numbers appeared in~Lemma~12. They are
uniquely determined by relation~(\ref{1.1}), where $\chi\in\mathcal
M'(\mathcal A,a,g)$ is arbitrarily chosen, taken together
with~(\ref{1.2}). This statement on uniqueness can actually be
strengthened as follows.\medskip

{\bf Lemma 12\,$'$.} {\it Given $\chi\in\mathcal M'(\mathcal
A,a,g)$, choose $\eta'_i\in\mathbb R$, $i\in I$, so that
$$
\sum_{i\in I}\alpha_i\eta'_i\geqslant\|\mathcal E(\mathcal
A,a,g)\|^2.
$$
If there holds\/}~(\ref{1.1}) {\it for $\eta'_i$ in place
of~$\eta_i$, then $\eta'_i=\eta_i$ for all $i\in I$.}\medskip

{\bf Proof.} This follows in the same manner as the uniqueness
statement in~Lemma~12.\hfill\phantom{B}\hfill$\Box$\medskip

{\bf 13.3.} The following assertion is specifying Lemma 12 for a
compact condenser~$\mathcal K$.\medskip

{\bf Corollary 16.} {\it Let $\mathcal A=\mathcal K$ be compact.
Given $\lambda_\mathcal K\in\mathcal S(\mathcal K,a,g)$, then for
every~$i$,}
\begin{align}
\alpha_ia_i\kappa(x,\lambda_\mathcal
K)&\geqslant\alpha_i\kappa(\lambda^i_\mathcal K,\lambda_\mathcal
K)\,g(x)&&\!\!\!\!\!\!\!\!\!\!\!\!\!\text{\it n.\,e. in \ } K_i,\label{hh}\\
\intertext{\it and hence} a_i\kappa(x,\lambda_\mathcal K)&
=\kappa(\lambda^i_\mathcal K,\lambda_\mathcal
K)\,g(x)&&\!\!\!\!\!\!\!\!\!\!\!\!\!\text{\it $\lambda^i_\mathcal
K$-almost everywhere.}\label{1.1''}
\end{align}

{\bf Proof.} In view of (\ref{S}) and (\ref{1.3'}), $\eta_i(\mathcal
K,a,g)$, $i\in I$, can be written in the form
$$\eta_i(\mathcal K,a,g)=\kappa(\lambda^i_\mathcal K,\lambda_\mathcal
K),$$ which leads to~(\ref{hh}) when substituted into~(\ref{1.1}).
Since $\lambda^i_\mathcal K$ has finite energy and is supported
by~$K_i$, the inequality in~(\ref{hh}) holds $\lambda^i_\mathcal
K$-almost everywhere in~$\mathbf X$. Hence, (\ref{1.1''}) must be
true, for if not, we would arrive at a contradiction by integrating
the inequality in~(\ref{hh}) with respect to~$\lambda^i_\mathcal
K$.\hfill\phantom{B}\hfill$\Box$\bigskip

{\bf\Large 14. Potentials of $\mathcal A$-vague cluster points of
minimizing nets}\nopagebreak\medskip

In this section we shall restrict ourselves to measures $\xi$ of the
class $\mathcal M_0(\mathcal A,a,g)$. It is clear from Corollary~7
that their potentials have all the properties described in Lemmas~12
and~12\,$'$. Our purpose is to show that, under proper additional
restrictions on the kernel, that description can be sharpened as
follows.\medskip

{\bf Lemma 13.} {\it In the case where $I^-\ne\varnothing$, assume
moreover that $\kappa(x,y)$ is continuous
on~$\overline{A^+}\times\overline{A^-}$, while\/}
$\kappa(\cdot,y)\to0$ ({\it as\/}~$y\to\infty$) {\it uniformly on
compact sets. Given $\xi\in\mathcal M_0(\mathcal A,a,g)$, then for
all $i\in I$,}
\begin{align} \alpha_ia_i\kappa(x,\xi)&
\geqslant\alpha_i\kappa(\xi^i,\xi)\,g(x)&&\!\!\!\!\!\!\!\!\!\!\!\!\!\!\!\!\!\!\text{\it
n.\,e. in \ }
A_i,\label{1.11}\\[4pt]
\alpha_i a_i\kappa(x,\xi)&\leqslant
\alpha_i\kappa(\xi^i,\xi)\,g(x)&&\!\!\!\!\!\!\!\!\!\!\!\!\!\!\!\!\!\!\text{\it
for all \ }
x\in S(\xi^i),\label{desc4'}\\
\intertext{\it and hence}
a_i\kappa(x,\xi)&=\kappa(\xi^i,\xi)\,g(x)&&\!\!\!\!\!\!\!\!\!\!\!\!\!\!\!\!\!\!\text{\it
n.\,e. in \ } A_i\cap S(\xi^i).\notag\end{align}

{\bf Proof.} Choose $\lambda_\mathcal K\in\mathcal S(\mathcal
K,a,g)$ such that $\xi$ is an $\mathcal A$-vague cluster point of
the net $(\lambda_\mathcal K)_{\mathcal K\in\{\mathcal K\}_\mathcal
A}$. Since this net belongs to $\mathbb M(\mathcal A,a,g)$, from
(\ref{1.3'}) and~(\ref{1.3''}) we get
$$
\eta_i=\kappa(\xi^i,\xi)=\lim_{\mathcal K\in\{\mathcal K\}_\mathcal
A}\,\kappa(\lambda^i_\mathcal K,\lambda_\mathcal K),\quad i\in I.
$$
Substituting this into~(\ref{1.1}) with~$\xi$ in place of~$\chi$
gives~(\ref{1.11}) as required.

We next proceed to prove~(\ref{desc4'}). To this end, fix $i$ (say
$i\in I^+$) and $x_0\in S(\xi^i)$. Without loss of generality it can
certainly be assumed that
\begin{equation}
\lambda_\mathcal K\to\xi\quad\mbox{$\mathcal A$-vaguely},
\label{l11}\end{equation} since otherwise we shall pass to a subnet
and change the notation. Then, due to~(\ref{1.1''}) and~(\ref{l11}),
one can choose $x_\mathcal K\in S(\lambda^i_\mathcal K)$ so that
\begin{equation}
x_\mathcal K\to x_0\quad\mbox{as \ } \mathcal K\uparrow\mathcal
A,\label{l112}
\end{equation}
$$
a_i\kappa(x_\mathcal K,\lambda_\mathcal K)=\kappa(\lambda^i_\mathcal
K,\lambda_\mathcal K)\,g(x_\mathcal K).
$$

Taking into account that, by~\cite[Lemma 2.2.1]{F1}, the map
$(x,\nu)\mapsto\kappa(x,\nu)$ is lower semicontinuous on the product
space $\mathbf X\times\mathfrak M^+$ (where $\mathfrak M^+$ is
equipped with the vague topology), we conclude from what has already
been shown that the desired relation~(\ref{desc4'}) will follow once
we prove
\begin{equation}
\kappa(x_0,R\xi^-)=\lim_{\mathcal K\in\{\mathcal K\}_\mathcal
A}\,\kappa(x_\mathcal K,R\lambda^-_\mathcal K).\label{hr}
\end{equation}

The case we are thus left with is $I^-\ne\varnothing$. Then,
according to our standing assumptions, $g_{\inf}>0$ and
$|a|<\infty$, and therefore there is $q\in(0,\infty)$ such that
\begin{equation}
R\lambda^-_\mathcal K(\mathbf X)\leqslant q\quad\mbox{for all \ }
\mathcal K\in\{\mathcal K\}_\mathcal A.\label{qq}
\end{equation}
Since, by (\ref{l11}) and Lemma 2, $R\lambda^-_\mathcal K\to R\xi^-$
vaguely, we thus get
\begin{equation}
R\xi^-(\mathbf X)\leqslant q.\label{xi}
\end{equation}

Fix $\varepsilon>0$. Under the assumptions of the lemma, one can
choose a compact neighborhood~$W_{x_0}$ of the point~$x_0$
in~$\overline{A^+}$ and a compact neighborhood~$F$ of~$W_{x_0}$
in~$\mathbf X$ so that
$$F_*:=F\cap\overline{A^-}\ne\varnothing$$
and
\begin{equation}
\bigl|\kappa(x,y)\bigr|<q^{-1}\varepsilon\quad\mbox{for all \ }
(x,y)\in W_{x_0}\times\complement F.\label{115}
\end{equation}

In the remainder, $\tilde{\complement}$ and~$\tilde{\partial}$
denote respectively the complement and the boundary of a set
relative to~$\overline{A^-}$ (where $\overline{A^-}$ is regarded to
be a topological subspace of~$\mathbf X$).

Having observed that $\kappa|_{W_{x_0}\times\overline{A^-}}$\, is
continuous, we proceed to construct a function
$$
\varphi\in\mathbf C_0(W_{x_0}\times\overline{A^-}\,)
$$
with the following properties:
\begin{equation}
\varphi|_{W_{x_0}\times F_*}=\kappa|_{W_{x_0}\times
F_*},\label{rest}
\end{equation}
\begin{equation}
\bigl|\varphi(x,y)\bigr|\leqslant q^{-1}\varepsilon\quad\mbox{for
all \ } (x,y)\in W_{x_0}\times\tilde{\complement}F_*.\label{118}
\end{equation}

To this end, consider a compact neighborhood~$V_*$ of~$F_*$
in~$\overline{A^-}$ and write
$$
f:=\left\{
\begin{array}{cl} \kappa & \mbox{on \ }
W_{x_0}\times\tilde{\partial}F_*,\\[2pt]
0 & \mbox{on \ } W_{x_0}\times\tilde{\partial}V_*.\\ \end{array}
\right.
$$
Note that $E:=(W_{x_0}\times\tilde{\partial}F_*)\cup
(W_{x_0}\times\tilde{\partial}V_*)$ is a compact subset of the
Hausdorff and compact, hence normal, space $W_{x_0}\times V_*$ and
$f$ is continuous on~$E$. By using the Tietze-Urysohn extension
theorem (see, e.\,g.,~\cite[Th.~0.2.13]{E2}), we deduce
from~(\ref{115}) that there exists a continuous function $\hat{f}: \
W_{x_0}\times V_*\to[-\varepsilon q^{-1},\varepsilon q^{-1}]$ such
that $\hat{f}|_E=f|_E$. Thus, the function in question can be
defined as follows:
$$
\varphi:=\left\{
\begin{array}{cl} \kappa & \mbox{on \ }
W_{x_0}\times F_*,\\[2pt]
\hat{f} & \mbox{on \ } W_{x_0}\times(V_*\setminus F_*),\\[2pt]
0 & \mbox{on \ } W_{x_0}\times\tilde{\complement}V_*.
\end{array}
\right.
$$

Furthermore, since the function $\varphi$  is continuous on
$W_{x_0}\times\overline{A^-}$ and has compact support, one can
choose a compact neighborhood~$U_{x_0}$ of~$x_0$ in~$W_{x_0}$ so
that
\begin{equation}
\bigl|\varphi(x,y)-\varphi(x_0,y)\bigr|<q^{-1}\varepsilon\quad\mbox{for
all \ } (x,y)\in U_{x_0}\times\overline{A^-}.\label{119}
\end{equation}

Given an arbitrary measure $\nu\in\mathfrak M^+(\,\overline{A^-}\,)$
with the property that $\nu(\mathbf X)\leqslant q$, we conclude
from~\mbox{(\ref{115})\,--\,(\ref{119})} that, for all $x\in
U_{x_0}$,
\begin{equation}
\bigl|\kappa\bigl(x,\nu|_{\complement
F}\bigr)\bigr|\leqslant\varepsilon,\label{121}
\end{equation}
\begin{equation}
\kappa\bigl(x,\nu|_{F}\bigr)=\int\varphi(x,y)\,d\bigl(\nu-\nu|_{\complement
F}\bigr)(y),\label{122}
\end{equation}
\begin{equation}
\Bigl|\int\varphi(x,y)\,d\nu|_{\complement
F}(y)\Bigr|\leqslant\varepsilon,\label{123}
\end{equation}
\begin{equation}
\Bigl|\int
\bigl[\varphi(x,y)-\varphi(x_0,y)\bigr]\,d\nu(y)\Bigr|\leqslant
\varepsilon.\label{124}
\end{equation}

Finally, choose $\mathcal K_0\in\{\mathcal K\}_\mathcal A$ so that
for all $\mathcal K\succ\mathcal K_0$ there hold $x_\mathcal K\in
U_{x_0}$ and
$$\Bigl|\int\varphi(x_0,y)\,d(R\lambda^-_\mathcal
K-R\xi^-)(y)\Bigr|<\varepsilon;$$ such a $\mathcal K_0$ exists
because of~(\ref{l11}) and~(\ref{l112}). Applying now relations
\mbox{(\ref{121})\,--\,(\ref{124})} to each of the
measures~$R\lambda^-_\mathcal K$ and~$R\xi^-$, which is possible due
to~(\ref{qq}) and~(\ref{xi}), for all $\mathcal K\succ\mathcal K_0$
we therefore get
\begin{equation*}
\begin{split}
\bigl|\kappa(x_\mathcal K&,R\lambda^-_\mathcal
K)-\kappa(x_0,R\xi^-)\bigr|\leqslant\bigl|\kappa\bigl(x_\mathcal
K,R\lambda^-_{\mathcal
K}\bigl|_{F}\bigr)-\kappa\bigl(x_0,R\xi^-\bigl|_{F}\bigr)\bigr|+2\varepsilon\\[7pt]
&{}\leqslant\Bigl|\int\varphi(x_\mathcal K,y)\,dR\lambda^-_\mathcal
K(y)-\int\varphi(x_0,y)\,dR\xi^-(y)\Bigr|+4\varepsilon\\[7pt]
&{}\leqslant\Bigl|\int\bigl[\varphi(x_\mathcal
K,y)-\varphi(x_0,y)\bigr]\,dR\lambda^-_\mathcal
K(y)\Bigr|+\Bigl|\int \varphi(x_0,y)\,d(R\lambda^-_\mathcal
K-R\xi^-)(y)\Bigr|+4\varepsilon\\[7pt]
&{}\leqslant\varepsilon+\varepsilon+4\varepsilon=6\varepsilon,
\end{split}
\end{equation*}
and (\ref{hr}) follows by letting $\varepsilon$ tend to~$0$. The
proof is complete.\hfill\phantom{B}\hfill$\Box$\bigskip

{\bf\Large 15. Proof of Theorems 2 and 3}\nopagebreak\medskip

We begin by showing that
\begin{equation} {\rm cap}\,(\mathcal A,a,g)\leqslant\|\hat{\Gamma}(\mathcal
A,a,g)\|^2.  \label{1} \end{equation} To this end,
$\|\hat{\Gamma}(\mathcal A,a,g)\|^2$ can certainly be assumed to be
finite. Then there are $\nu\in\hat{\Gamma}(\mathcal A,a,g)$ and
$\mu\in\mathcal E^0(\mathcal A,a,g)$, the existence of~$\mu$ being
clear from~(\ref{nonzero1}) and Corollary~3.
By~\cite[Lemma~2.3.1]{F1}, the inequality in~(\ref{adm1}) holds
$\mu^i$-almost everywhere. Integrating it with respect to~$\mu^i$
and then summing up over all $i\in I$, in view of $\int
g\,d\mu^i=a_i$ we get
$$ \kappa(\nu,\mu)\geqslant\sum_{i\in I}\,c_i(\nu), $$ hence
$\kappa(\nu,\mu)\geqslant1$ by~(\ref{adm2}), and finally $$
\|\nu\|^2\|\mu\|^2\geqslant1
$$
by the Cauchy-Schwarz inequality. The last relation, being valid for
arbitrary $\nu\in\hat{\Gamma}(\mathcal A,a,g)$ and $\mu\in\mathcal
E^0(\mathcal A,a,g)$, implies~(\ref{1}), which in turn yields
Theorem~2 provided ${\rm cap}\,\mathcal A=\infty$.

We are thus left with proving both Theorems~2 and~3 in the case
where ${\rm cap}\,\mathcal A$ is finite. Then the $\mathcal
E(\mathcal A,a\,{\rm cap}\,\mathcal A,g)$-problem can be considered
as well.

Taking (\ref{def}) and (\ref{iden}) into account, we deduce from
Lemmas~8 and~12 with~$a$ replaced by $a\,{\rm cap}\,\mathcal A$
that, for every $\chi\in\mathcal M'(\mathcal A,a\,{\rm
cap}\,\mathcal A,g)$,
\begin{equation}
\|\chi\|^2={\rm cap}\,\mathcal A \label{2}
\end{equation}
and there exist unique $\tilde{\eta}_i\in\mathbb R$, $i\in I$, such
that
\begin{equation}
\alpha_ia_i\kappa(x,\chi)\geqslant\tilde{\eta_i}\,g(x)\quad\mbox{n.\,e.~in
\ } A_i,\quad i\in I,\label{tilde1}
\end{equation}
\begin{equation}
\sum_{i\in I}\,\tilde{\eta_i}=1.\label{tilde2}
\end{equation}
Actually,
\begin{equation}
\tilde{\eta_i}=\alpha_i\,{\rm cap}\,\mathcal A^{-1}\,\eta_i(\mathcal
A,a\,{\rm cap}\,\mathcal A,g),\quad i\in I,\label{tilde3}
\end{equation}
where $\eta_i(\mathcal A,a\,{\rm cap}\,\mathcal A,g)$, $i\in I$, are
the numbers uniquely determined in~Sec.~13.

Using the property of sets of interior capacity zero mentioned
in~Sec.~7.1 and the fact that the potentials of equivalent
in~$\mathcal E$ measures coincide nearly everywhere in~$\mathbf X$,
we conclude from~(\ref{tilde1}) and~(\ref{tilde2}) that
$$
\mathcal M'_\mathcal E(\mathcal A,a\,{\rm cap}\,\mathcal
A,g)\subset\hat{\Gamma}(\mathcal A,a,g).
$$
Together with~(\ref{1}) and~(\ref{2}), this implies that, for every
$\sigma\in\mathcal M'_\mathcal E(\mathcal A,a\,{\rm cap}\,\mathcal
A,g)$,
$${\rm cap}\,\mathcal A=\|\sigma\|^2\geqslant\|\hat{\Gamma}(\mathcal
A,a,g)\|^2\geqslant{\rm cap}\,\mathcal A,$$ which completes the
proof of Theorem~2. The last two relations also yield
$$
\mathcal M'_\mathcal E(\mathcal A,a\,{\rm cap}\,\mathcal
A,g)\subset\hat{\mathcal G}(\mathcal A,a,g).
$$
As both the sides of this inclusion are equivalence classes
in~$\mathcal E$ (see Lemma~9), they must actually be equal, and
(\ref{desc}) follows.

Applying Lemma~12\,$'$ for~$a\,{\rm cap}\,\mathcal A$ in place
of~$a$, we deduce from~(\ref{desc}) that $c_i(\hat{\omega})$, $i\in
I$, satisfying (\ref{adm1}) and~(\ref{adm2}) for
$\nu=\hat{\omega}\in\hat{\mathcal G}(\mathcal A,a,g)$, are
determined uniquely, do not depend on the choice of~$\hat{\omega}$,
and are actually equal to~$\tilde{\eta_i}$. Therefore,
substituting~(\ref{1.3'}) and, subsequently, (\ref{1.3''})
for~$a\,{\rm cap}\,\mathcal A$ in place of~$a$ into~(\ref{tilde3}),
we get~(\ref{snt1}) and (\ref{cnt}). This proves
Theorem~3.\hfill\phantom{B}\hfill$\Box$\bigskip

{\bf\Large 16. Proof of Theorem 5}\nopagebreak\medskip

We start by observing that $\mathcal D(\mathcal A,a,g)$ is nonempty,
contained in an equivalence class in~$\mathcal
E(\,\overline{\mathcal A}\,)$, and satisfies the inclusions
\begin{equation}
\mathcal D(\mathcal A,a,g)\subset\mathcal M(\mathcal A,a\,{\rm
cap}\,\mathcal A,g)\subset\mathcal M'(\mathcal A,a\,{\rm
cap}\,\mathcal A,g)\cap\mathcal E(\,\overline{\mathcal
A},\leqslant\!a\,{\rm cap}\,\mathcal A,g). \label{DD}
\end{equation}
Indeed, this follows from (\ref{D}), Corollary~7, and Lemma~8, the
last two being taken for~$a\,{\rm cap}\,\mathcal A$ in place of~$a$.
Substituting~(\ref{desc2}) into~(\ref{DD}) gives~(\ref{gammain}) as
required.

Since, by (\ref{gammain}), every $\gamma\in\mathcal D(\mathcal
A,a,g)$ is a minimizer in the $\Gamma(\mathcal A,a,g)$-problem, the
claimed relations~(\ref{5.1}) and~(\ref{5.3}) are obtained directly
from Theorem~3 and~4 in view of Definition~6. To show that
$C_i(\mathcal A,a,g)$, $i\in I$, can actually be given by means
of~(\ref{5.2}), one only needs to substitute~$\gamma$ instead
of~$\zeta$ into~(\ref{snt1})~--- which is possible due
to~(\ref{DD})~--- and use Corollary~8.

Assume for a moment that, if $I^-\ne\varnothing$, then the kernel
$\kappa(x,y)$ is continuous on~$\overline{A^+}\times\overline{A^-}$,
while $\kappa(\,\cdot\,,y)\to0$ (as $y\to\infty$) uniformly on
compact sets. In order to establish~(\ref{5.4}), it suffices to
apply Lemma~13 (with~$a\,{\rm cap}\,\mathcal A$ in place of~$a$)
to~$\gamma$, which can be done because of~(\ref{D}), and then
substitute~(\ref{5.2}) into the result obtained.

To prove that $\mathcal D(\mathcal A,a,g)$ is $\mathcal A$-vaguely
compact, fix $(\gamma_s)_{s\in S}\subset\mathcal D(\mathcal A,a,g)$.
Then the inclusion~(\ref{gammain}) and Lemma~10 yield that this net
is $\mathcal A$-vaguely bounded and hence, by~Lemma~1, $\mathcal
A$-vaguely relatively compact. Let $\gamma_0$ denote one of its
$\mathcal A$-vague cluster points, and let $(\gamma_t)_{t\in T}$ be
a subnet of~$(\gamma_s)_{s\in S}$ that converges $\mathcal
A$-vaguely to~$\gamma_0$. In view of~(\ref{D}), the proof will be
completed once we show that
\begin{equation}
\gamma_0\in\mathcal M_0(\mathcal A,a\,{\rm cap}\,\mathcal
A,g).\label{GG}
\end{equation}

By~(\ref{D}), for every $t\in T$ there exist a subnet $(\mathcal
K_{s_t})_{s_t\in S_t}$ of the net $(\mathcal K)_{\mathcal
K\in\{\mathcal K\}_\mathcal A}$ and
$$\lambda_{s_t}\in\mathcal
S(\mathcal K_{s_t},a\,{\rm cap}\,\mathcal A,g),\quad s_t\in S_t,$$
such that $\lambda_{s_t}$ approaches $\gamma_t$ $\mathcal A$-vaguely
as $s_t$ ranges over~$S_t$. Consider the Cartesian product
$\prod\,\{S_t: t\in T\}$~--- that is, the collection of all
functions~$\psi$ on~$T$ with $\psi(t)\in S_t$, and let~$D$ denote
the directed product $T\times\prod\,\{S_t: t\in T\}$ (see,
e.\,g.,~\cite[Chap.~2,~\S~3]{K}). Given $(t,\psi)\in D$, write
$$
\mathcal K_{(t,\psi)}:=\mathcal K_{\psi(t)}\quad\mbox{and}\quad
\lambda_{(t,\psi)}:=\lambda_{\psi(t)}.
$$
Then the theorem on iterated limits from \cite[Chap.~2, \S~4]{K}
yields that $(\lambda_{(t,\psi)})_{(t,\psi)\in D}$ converges
$\mathcal A$-vaguely to~$\gamma_0$. Since, as can be seen from the
above construction, $(\mathcal K_{(t,\psi)})_{(t,\psi)\in D}$ forms
a subnet of $(\mathcal K)_{\mathcal K\in\{\mathcal K\}_\mathcal A}$,
this proves~(\ref{GG}) as
required.\hfill\phantom{B}\hfill$\Box$\bigskip

{\bf\Large 17. Proof of Proposition 2}\nopagebreak\medskip

Consider $\nu\in\mathcal E(\,\overline{\mathcal A}\,)$ and
$\tau_i\in\mathbb R$, $i\in I$, satisfying both the
assumptions~(\ref{desc3'}) and~(\ref{adm2''}), and fix arbitrarily
$\gamma_\mathcal A\in\mathcal D(\mathcal A,a,g)$ and $(\mu_t)_{t\in
T}\in\mathbb M(\mathcal A,a\,{\rm cap}\,\mathcal A,g)$. Since
$\mu^i_t$ is concentrated on~$A_i$ and has finite energy and compact
support, the inequality in~(\ref{desc3'}) holds $\mu^i_t$-almost
everywhere. Integrating it with respect to~$\mu^i_t$ and then
summing up over all $i\in I$, in view of~(\ref{5.1})
and~(\ref{adm2''}) we obtain
$$
2\,\kappa(\mu_t,\nu)\geqslant\|\gamma_\mathcal A\|^2+\|\nu\|^2,\quad
t\in T.
$$
But $(\mu_t)_t\in T$ converges to~$\gamma_\mathcal A$ in the strong
topology of the semimetric space~$\mathcal E(\,\overline{\mathcal
A}\,)$, which is clear from~(\ref{DD}) and Lemma~8 with~$a\,{\rm
cap}\,\mathcal A$ instead of~$a$. Therefore, passing in the
preceding relation to the limit as $t$ ranges over~$T$, we get
$$
\|\nu-\gamma_\mathcal A\|^2=0,
$$
which is a part of the conclusion of the proposition. In turn, the
preceding relation implies that, actually, the right-hand side
in~(\ref{adm2''}) is equal to~$1$, and that $\nu\in\mathcal
M'(\mathcal A,a\,{\rm cap}\,\mathcal A,g)$. Since, in view
of~Theorem~3, the latter means that
$$
R\nu\in\hat{\mathcal G}(\mathcal A,a,g),
$$
the claimed relation (\ref{un})
follows.\hfill\phantom{B}\hfill$\Box$\bigskip

{\bf\Large 18. Proof of Theorem 6}\nopagebreak\medskip

To establish~(\ref{6.1}), fix $\mu\in\mathcal E(\mathcal A,a,g)$.
Under the assumptions of the theorem, either $g_{\inf}>0$, and
consequently $\mu^i(\mathbf X)<\infty$ for all $i\in I$, or $\mathbf
X$ is countable at infinity; in any case, every~$A_i$, $i\in I$, is
contained in a countable union of $\mu^i$-integrable sets.
Therefore, by~\cite{B2,E2} (cf.~also the appendix below),
\begin{align*}
\int g\,d\mu^i&=\lim_{n\in\mathbb N}\,\int
g\,d\mu^{i}_{\mathcal A_n},&&\!\!\!\!\!\!\!\!\!\!\!\!\!\!\!\!\!\!\!\!\!\!\!\!\!\!\!\!\!\!i\in I,\\
\kappa(\mu^i,\mu^j)&=\lim_{n\in\mathbb N}\,\kappa(\mu^{i}_{\mathcal
A_n},\mu^{j}_{\mathcal
A_n}),&&\!\!\!\!\!\!\!\!\!\!\!\!\!\!\!\!\!\!\!\!\!\!\!\!\!\!\!\!\!\!i,\,j\in
I,
\end{align*}
where $\mu^{i}_{\mathcal A_n}$ denotes the trace of $\mu^{i}$
upon~$A^i_n$. Applying the same arguments as in the proof of
Lemma~5, but now based on the preceding two relations instead
of~(\ref{w}) and~(\ref{ww}), we arrive at~(\ref{6.1}) as required.

By~(\ref{nonzero1}) and~(\ref{6.1}), for every $n\in\mathbb N$,
${\rm cap}\,(\mathcal A_n,a,g)$ can certainly be assumed to be
nonzero. Suppose moreover that ${\rm cap}\,(\mathcal A,a,g)$ is
finite; then, by~(\ref{increas'}), so is ${\rm cap}\,(\mathcal
A_n,a,g)$. Hence, according to Theorem~5, there exists
\begin{equation}
\gamma_n:=\gamma_{\mathcal A_n}\in\mathcal D(\mathcal
A_n,a,g).\label{17.1}\end{equation}

Observe that $R\gamma_n$ is a minimizer in the
$\hat{\Gamma}(\mathcal A_n,a,g)$-problem, which is clear
from~(\ref{desc}), (\ref{desc2}), and~(\ref{gammain}). Since,
furthemore,
$$\hat{\Gamma}(\mathcal A_{n+1},a,g)\subset\hat{\Gamma}(\mathcal A_n,a,g),$$
application of Lemma~4 to $\mathcal H=\hat{\Gamma}(\mathcal
A_n,a,g)$, $\nu=R\gamma_{n+1}$, and $\lambda=R\gamma_n$ gives
$$
\|\gamma_{n+1}-\gamma_n\|^2\leqslant\|\gamma_{n+1}\|^2-\|\gamma_n\|^2.
$$
Also note that $\|\gamma_n\|^2$, $n\in\mathbb N$, is a Cauchy
sequence in~$\mathbb R$, because, by~(\ref{6.1}), its limit exists
and, being equal to~${\rm cap}\,\mathcal A$, is finite. The
preceding inequality therefore yields that $(\gamma_n)_{n\in\mathbb
N}$ is a strong Cauchy sequence in the semimetric space~$\mathcal
E(\,\overline{\mathcal A}\,)$.

Besides, since ${\rm cap}\,\mathcal A_n\leqslant{\rm cap}\,\mathcal
A$, we derive from~(\ref{gammain}) that
$$(\gamma_n)_{n\in\mathbb N}\subset\mathcal E(\,\overline{\mathcal
A}\,,\leqslant\!a\,{\rm cap}\,\mathcal A,g).$$ Hence, by Theorem~7,
there exists an $\mathcal A$-vague cluster point~$\gamma$
of~$(\gamma_n)_{n\in\mathbb N}$, and
$$
\lim_{n\in\mathbb N}\,\|\gamma_n-\gamma\|^2=0.
$$
Let $(\gamma_t)_{t\in T}$ denote a subnet of the sequence
$(\gamma_n)_{n\in\mathbb N}$ that converges $\mathcal A$-vaguely and
strongly to~$\gamma$. We next proceed to show that
\begin{equation}
\gamma\in\mathcal D(\mathcal A,a,g). \label{theta}\end{equation}

For every $t\in T$, consider the ordered family $\{\mathcal
K_t\}_{\mathcal A_t}$ of all compact condensers $\mathcal
K_t\prec\mathcal A_t$. By~(\ref{17.1}), there exist a subnet
$(\mathcal K_{s_t})_{s_t\in S_t}$ of $(\mathcal K_t)_{\mathcal
K_t\in\{\mathcal K_t\}_{\mathcal A_t}}$ and
$$\lambda_{s_t}\in\mathcal S(\mathcal K_{s_t},a\,{\rm
cap}\,\mathcal K_{s_t},g)$$ such that $(\lambda_{s_t})_{s_t\in S_t}$
converges both strongly and $\mathcal A$-vaguely to~$\gamma_t$.
Consider the Cartesian product $\prod\,\{S_t: t\in T\}$, that is,
the collection of all functions~$\psi$ on~$T$ with $\psi(t)\in S_t$,
and let $D$ denote the directed product $T\times\prod\,\{S_t: t\in
T\}$. Given $(t,\psi)\in D$, write
$$
\mathcal K_{(t,\psi)}:=\mathcal K_{\psi(t)}\quad\mbox{and}\quad
\lambda_{(t,\psi)}:=\lambda_{\psi(t)}.
$$
Then the theorem on iterated limits from~\cite[Chap.~2, \S~4]{K}
yields that $(\lambda_{(t,\psi)})_{(t,\psi)\in D}$ converges both
strongly and $\mathcal A$-vaguely to~$\gamma$. Since $(\mathcal
K_{(t,\psi)})_{(t,\psi)\in D}$ forms a subnet of~$(\mathcal
K)_{\mathcal K\in\{\mathcal K\}_\mathcal A}$, this
proves~(\ref{theta}) as required.

What is finally left is to prove~(\ref{6.3}). By Corollary~13, for
every $n\in\mathbb N$ one can choose a compact condenser $\mathcal
K^0_n\prec\mathcal A_n$ so that
$$
\bigl|C_i(\mathcal A_n,a,g)-C_i(\mathcal
K^0_n,a,g)\bigr|<n^{-1},\quad i\in I.
$$
This $\mathcal K^0_n$ can certainly be chosen so large that the
sequence obtained, $(\mathcal K^0_n)_{n\in\mathbb N}$, forms a
subnet of $(\mathcal K)_{\mathcal K\in\{\mathcal K\}_\mathcal A}$;
therefore, repeated application of Corollary~13 yields
$$\lim_{n\in\mathbb N}\,C_i(\mathcal
K^0_n,a,g)=C_i(\mathcal A,a,g).$$ This leads to (\ref{6.3}) when
combined with the preceding
relation.\hfill\phantom{B}\hfill$\Box$\bigskip

{\bf\Large 19. Acknowledgments}\nopagebreak\medskip

The author is greatly indebted to P.~Dragnev, D.~Hardin, and
E.~B.~Saff for several valuable comments concerning this
study.\bigskip

{\bf\Large 20. Appendix}\nopagebreak\medskip

Let $\nu\in\mathfrak M^+(\mathbf X)$ be given. As in~\cite[Chap.~4,
\S~4.7]{E2}, a set~$E\subset\mathbf X$ is called
\mbox{$\nu$-$\sigma$}-{\it fi\-ni\-te\/} if it can be written as a
countable union of $\nu$-integr\-able sets.

The following assertion, related to the theory of measures and
integration, has been used in Sec.~18. Although it is not difficult
to deduce it from~\cite{B2,E2}, we could not find there a proper
reference.\medskip

{\bf Lemma 14.} {\it Consider a lower semicontinuous function~$\psi$
on~$\mathbf X$ such that $\psi\geqslant0$ unless the space~$\mathbf
X$ is compact, and let $E$ be the union of an increasing sequence of
$\nu$-measur\-able sets~$E_n$, $n\in\mathbb N$. If moreover $E$ is
\mbox{$\nu$-$\sigma$}-fi\-ni\-te, then}
$$
\int\psi\,d\nu_E=\lim_{n\in\mathbb N}\,\int\psi\,d\nu_{E_n}.
$$

{\bf Proof.} Without loss of generality, we can certainly assume
$\psi$ to be nonnegative. Then for every
\mbox{$\nu$-$\sigma$}-fi\-ni\-te set~$Q$,
\begin{equation}
\int\psi\,d\nu_Q=\int\psi\varphi_Q\,d\nu, \label{Q}
\end{equation}
where $\varphi_Q(x)$ equals~$1$ if $x\in Q$, and~$0$ otherwise.
Indeed, this can be concluded from~\cite[Chap.~4, \S~4.14]{E2} (see
Propositions~4.14.1 and~4.14.6).

On the other hand, since $\psi\varphi_{E_n}$, $n\in\mathbb N$, are
nonnegative and form an increasing sequence with the upper
envelope~$\psi\varphi_E$, \cite[Prop.~4.5.1]{E2} gives
$$
\int\psi\varphi_E\,d\nu=\lim_{n\in\mathbb
N}\,\int\psi\varphi_{E_n}\,d\nu.
$$
Applying (\ref{Q}) to both the sides of this equality, we obtain the
lemma.\hfill\phantom{B}\hfill$\Box$

\bigskip
\flushleft

{\small Institute of Mathematics\\
National Academy of Sciences of Ukraine\\
3 Tereshchenkivska Str.\\
01601, Kyiv-4, Ukraine\\
e-mail: natalia.zorii@gmail.com}

\end{document}